\documentclass[11pt]{article}
\usepackage{lscape}
\renewcommand{\baselinestretch}{1.5}
\title{Computation of Irreducible Characters of $GL(5,q)$}
\author{By Shiv K. Gupta}
\date{\mbox{}}
\begin{document}
\maketitle
%
%
\section{INTRODUCTION}\footnote{{\em Mathematics Subject Classification.}
Primary 20C15; Secondary 20C30}
In this paper we determine the ordinary irreducible characters of the
5-dimensional full linear group over a Galois field of $q$ elements. We
use the techniques developed by Green [1].

The characters of $GL(2,q)$ have been known for long time. Those of
$GL(3,q)$
and $GL(4,q)$ have been determined by Steinberg [11]. Steinberg has also
developed some techniques to find some of the characters in the general
case.

There are 42 ``types'' of characters of $GL(5,q)$. Each of these types
contains
a number of characters which is a function of $q$. The total number of
characters of $GL(5,q)$ is $q^5-q^2-q+1$. There is a perfect symmetry in the
description of conjugate classes and characters.

Thee 42 ``types'' of characters may in turn be grouped into 17 classes, such
that the characters belonging to the same class involve the same basic
constituents--the basic characters. For example, the characters of the first
seven types $A'_1$ to $A'_7$ (each containing $q-1$ characters) form one class.
Each of these types involve the same basic constituents--the basic characters
of type $A'$ , occurring with distinct coefficients. [These coefficients being
the linear combinations which involve what are known as Green Polynomials,
with (dual) functions which involve the characters of the symmetric groups.]
%
\newcommand{\ea}{\mbox{$\epsilon_{1}$}}
\newcommand{\eb}{\mbox{$\epsilon_{2}$}}
\newcommand{\ec}{\mbox{$\epsilon_{3}$}}
\newcommand{\ed}{\mbox{$\epsilon_{4}$}}
\newcommand{\ee}{\mbox{$\epsilon_{5}$}}

\section{Definitions and Notation}
We have used the standard notation throughout. $GF(q)$ will denote a finite
field of $q$ elements; $GL(n,q)$ the group of $n\times n$ non-singular
matrices over GF($q$), $S_{n}$ the symmetric group on $n$  letters.
If $\lambda$ is a
partition of $n$ then we shall write $\mid\lambda\mid=n$. Conjugate class
``types'' of $GL(5,q)$ have been denoted by the letters $A,B,C,\ldots ,Q$
and the characters types by $A', B',C',\ldots , Q'$.

Let $\rho = \{\ell_{1},\ell_{2},\ldots\} , \ell_{1} \geq \ell_{2}\cdots $
be a partition
of n. Then the partition $\rho^{'} = \{\ell_{1}^{'}, \ell_{2}^{'}, \ldots\}$ is
said to be conjugate to the partition $\rho$ if
$ \ell_{r}^{'}$ = number of parts of $\rho$ which are $\geq r$.
A conjugate class $c$ of $GL(n,q)$ is said to be {\em primary} if its
characteristic polynomial is a power of irreducible (over $GF(q)$)
polynomial. We shall number the seven partitions for $n = 7$ as follows:

\begin{displaymath}
   \rho_{1} = 1^{5},\  \rho_{2}=2.1^{3},\  \rho_{3} = 2^{2}.1,\
    \rho_{4}=3.1^{2},\  \rho_{5} = 3.2,\  \rho_{6}=4.1,\  \rho_{7} = 5.
\end{displaymath}
$\omega_{i} , 1\leq i\leq 5$ will denote the primitive elements of GF($q^{i}$)
such that
\begin{displaymath}
   \omega_{1} = \omega_{2}^{q+1} = \omega_{3}^{q^{2}+q+1} =
   \omega_{4}^{q^{3}+q^{2}+q+1} = \omega_{5}^{q^{4}+q^{3}+q^{2}+q+1}.
\end{displaymath}
$\varepsilon_{i}, 1\leq i\leq 5$ are the primitive ($q^{i}$-1)-$th$
roots of unity such that
\begin{displaymath}
   \varepsilon_{1} =
    \varepsilon_{2}^{q+1} =
    \varepsilon_{3}^{q^{2}+q+1} =
    \varepsilon_{4}^{q^{3}+q^{2}+q+1} =
    \varepsilon_{5}^{q^{4}+q^{3}+q^{2}+q+1}
\end{displaymath}
Let $f$ be an irreducible polynomial of degree $s$ over GF($q$). Let $\omega$
be a generating element of GF($q^{s}$). Then a root of $f$ has the form
$\omega^{k}$ where $k$ is uniquely determined modulo $(q^{s}-1)$, and the set
$\{\omega^{k}, \omega^{kq}, \omega^{kq^{2}},\cdots,\omega^{kq^{s-1}}\}$
consists of all the roots of $f$ since $k, kq, kq^{2}, \ldots, kq^{s-1}$ are
distinct residue mod $(q^{s}-1)$. Following Green [1] we shall call the
set $\{k,kq,kq^{2},\ldots,kq^{s-1}\}$ an $s$-$simplex$ with
$k,kq,kq^{2},\ldots,kq^{s-1}$ as its roots and $s$ its degree.
$\chi_{\rho}^\lambda$ will denote the value of the irreducible character of
$S_{\mid\lambda\mid}$ corresponding to the partition $\lambda$ of
$\mid\lambda\mid$ , on the class corresponding to the partition $\rho$.
$\sum_{i,j,k\ldots}$ will mean that the summation is taken over all the
permutations of $i,j,k\ldots$.

Let $G$ denote the set of all simplexes of degrees $\leq n$, and $P$ be the
set of all partitions of numbers $\leq n$ . Then a function
\mbox{$\nu : G\rightarrow P$} such that
\begin{displaymath}
   \sum_{g\in G}\mid\nu(g)\mid deg(g) = n
\end{displaymath}
specifies what we shall call a {\em dual class} and will designate it by the
symbol $e = (\ldots g^{\nu(g)}\ldots)$. To each symbol of a dual class
there will correspond a character which we shall denote by the same
symbol.

\section{Description of Conjugate Classes}

Let $f(t) = t^{d}-a_{d-1}-\cdots-a_{0}$ be a polynomial over $GF(q)$.
Following Green [1] we define
\begin{displaymath}
 U(f) = \left(
  \begin{array}{cccccccc}
                               & 1 & & & & & &\\
                               & & . & & & & &\\
                               & & & . & & & &\\
                               & & & & . & & &\\
                               & & & & & . & &\\
                               & & & & & & 1 &\\
                               a_{0} & a_{1} & a_{2} & . & . & . & . & a_{d-1}
  \end{array}
     \right)
\end{displaymath}
and

\begin{displaymath}
U_{m}(f) = \left(
  \begin{array}{cccccccc}
                            U(f) & I_{d} & & & & & \\
                               & U(f) & I_{d}   & & & &\\
                               & & & . & & & &\\
                               & & & & & .  & &\\
                               & & & & & & .  &\\
                               & & & & & & & I_{d}\\
                               & & & & & & & U(f)
  \end{array}
     \right)
\end{displaymath}
where $I_{d}$ is an $d\times d$ identity matrix and there are $m$ blocks of
$U(f)$.

Let $\lambda = \{\ell_{1}, \ell_{2},\ldots,\ell_{p} \} \
\ell_{1} \geq \ell_{2} \ldots \geq \ell_{p}$ , be a partition of $n$.
We define
\[
U_{\lambda}(f) =
diag\{U_{\ell_{1}}(f),U_{\ell_{2}}(f),\ldots,U_{\ell_{r}}(f)\}.
\]
Then
$U_{\lambda}(f)$ is a matrix with characteristic polynomial
$f(t)^{\mid\lambda\mid}$. Let $A \in$ GL($n,q$) have characteristic
polynomial $f_{1}^{n_{1}}.f_{2}^{n_{2}}\cdots f_{k}^{n_{k}}$ where
$f_{1}, f_{2},\cdots,f_{k}$ are distinct irreducibles polynomials
over GF($q$). Then $\sum_{i=1}^{k} n_{i} deg(f_{i}) = n $, where
deg($f_{i}$) is the degree of the polynomial $f_{i}$. The matrix
$A\in GL(n,q)$ is conjugate to a matrix
diag\{$U_{\nu_{1}}(f_{1}),U_{\nu_{2}}(f_{2}),\cdots,U_{\nu_{k}}(f_{k})$\},
where $\nu_{1},\nu_{2},\cdots,\nu_{k}$ are certain partitions of
$n_{1}, n_{2},\cdots,n_{k}$ respectively.

Conversely, let $F$ denote the set of all irreducible polynomials over
GF($q$) of degree $\leq n$.  A conjugate class of GL($n,q$) can be
specified by a partition valued function \mbox{$\nu : F\rightarrow P$}
such that \mbox{$\sum_{f \in F}\mid \nu(f)\mid deg(f) = n$}. We designate
this class by the symbol \mbox{$(\cdots f^{\nu(f)} \cdots)$}.

In the case of $GL(5,q)$ we are thus able to describe the conjugate classes
systematically. There are $q^5-q^2-q+1$ conjugate classes in $GL(5,q)$. These
are divided into 42  ``types'' and all these  types do appear  if $q\geq 7$.
Each type consists of a number of conjugate classes and this number is a
function of $q$. All the classes of a particular type have similar minimal
polynomial. In $GL(5,2), GL(5,3) , GL(5,4)$ and $GL(5,5)$ there are
respectively
17, 32, 39 and 41 types of classes. We give a list of the canonical forms of
the various types of conjugate classes in table 1. The table 2 gives the
number of classes in each type and also the number of elements in each class.

\section{Uniform Functions and Basic Uniform Functions}

Let $A$ denote a square matrix of order $n$ with coefficients in GF($q$).
Let $V = V_n$ denote an n-dimensional vector space over GF($q$).
$V$ can be considered as a GF($q$)[X]-module denoted by $V_A$ by defining:

\begin{displaymath}
(a_0 + a_1 X + a_2 X^2 + \cdots + a_r X^r) v =
a_0+ a_1 Av + a_2 A^2v+\cdots+a_rA^r v
\end{displaymath}
where $a_i \in$ GF$(q) , v \in V$.
Let $\lambda$ be a fixed partition of $n , f(t) = t$ and $A=U_{\lambda}(f)$.
Set $V_{\lambda}=V_{U_{\lambda}}(f)$. If the partition $\lambda$ has $p$
parts then we define
%
\[
   k(\lambda ,q) = \left\{ \begin{array}{ll}
        (1-q)(1-q^2)(1-q^3) \cdots (1-q^{p-1})  & \mbox{if $p > 1$} \\
         0                                     & \mbox{if $p = 1$}
                          \end{array}
                   \right.
\]
%
Let $\lambda_1, \lambda_2,\cdots,\lambda_k$ be any partitions such that
\mbox {$\sum_{i=1}^{i=k} \mid \lambda_{i} \mid = \mid \lambda \mid $}
Then \mbox{$ g^{\lambda}_{\lambda_1,\lambda_2,\ldots,\lambda_k}$}
is the number of chains
\mbox{$  V_{\lambda} =
                   V^0   \supset
                   V^1   \supset
                   V^2   \supset
                   \cdots\supset
                   V^{k-1}\supset
                   V^{k} = 0
     $}
of sub-modules of $V_{\lambda}$ in which the factor $V^{(i-1)}/V^i$ is
isomorphic to $V_{\lambda_{i}}, i = 1,2,\cdots,k.$
\subsection{Green Polynomial}

Let $\lambda, \rho$ be partitions of $n$. Assume that
$\rho= 1^{r_1}.2^{r_2}.3^{r_3}\cdots$  \
the polynomials
\[
Q^{\lambda}_{\rho}(q) = \sum g^{\lambda}_{\lambda_1,\lambda_2\cdots}(q)
.k(\lambda_1,q).k(\lambda_2,q)\cdots
\]
the sum being taken over all sequences ($\lambda_1,\lambda_2,\cdots$) of
partitions of $\rho$ such that
$\lambda_1, \lambda_2,\cdots,\lambda_{r_1}$ are partitions of 1 and
$\lambda_{r_1+1},\lambda_{r_1+2},\cdots,\lambda_{r_2}$ are
partitions of $2$ etc.
are known as {\em Green Polynomials}.
The Green Polynomials are closely related with the characters of the
symmmetric groups. In fact it is known [3] that if
$\overline{\lambda}=(\overline{\lambda_1},\overline{\lambda_2},\cdots)$
is a partition conjugate to $\lambda$ and
$n_{\lambda}= \frac{1}{2}\sum\overline{\lambda_i}(\overline{\lambda_i}$-1)
then
$\chi^{\lambda}_{\rho}$ is equal to the coefficient of $q^{n_{\lambda}}$ in
$Q^{\lambda}_{\rho}$. Green Polynomials satisfy the following orthogonality
relations:
\begin{enumerate}
   \item If $\nu$ be any partition of $n$ then
      \[
         Q^{\nu}_{\rho+\sigma}(q) =
             \sum g^{\nu}_{\lambda,\mu}(q)
             Q^{\lambda}_{\rho}(q) Q^{\mu}_{\sigma}(q)
      \]
       where the summation is over all pairs $\lambda, \mu$ of partitions
       $\ell$ and $m$, where $\ell+m = n$. Here if
       $\rho=1^{r_1}.2^{r_2}.3^{r_3}\cdots$ and
       $\sigma=1^{s_1}.2^{s_2}.3^{s_3}\cdots$ then the partition
       $\rho+\sigma$ is defined to be the partition
       $1^{r_1+s_1}.2^{r_2+s_2}.3^{r_3+s_3}\cdots$
   \item Let $\nu , \rho$ be partitions of $n$  then
       \[
          \sum_{\mid\rho\mid=n}\frac{1}{z}Q^{\nu}_{\rho}(q)=1
       \]
    where $z_{\rho} = 1^{r_1}(r_{1}!.2^{r_2}(r_{2}!)\cdots$ is the order
    of the centralizer (in $S_n$) of an element of an element in $S_n$
    corresponding to the partition $\rho$.
\end{enumerate}

Green [1] gives a table of the polynomials $Q^{\lambda}_{\rho}$ for
$1\leq n\leq 5$. Morris's paper [7] deals mainly with these polynomials
and a method of their determination. He also gives a table of
these polynomials for $n=6,7$.

\subsection{$\rho$-Variables and $\rho$-Functions}

Let $\rho= 1^{r_1}.2^{r_2}\cdots$ be a partition of $\rho$. Then a
formal set $X^{\rho}= \{x_{11},\cdots,x_{1r_1},x_{21},\cdots,x_{2r_2}\cdots\}$
of variables is called a {\em set o}f $\rho$-{\em variables}, and each member
of the set is called a $\rho$-{\em variable}. For each positive integer $d$
there are $r_d\ \rho$-variables namely $x_{d1},x_{d2},\cdots,x_{dr_d}$.
Each of these $r_d\ \rho$-variables is said to have degree deg($x_{di})=d$,
$i=1,2,\cdots,r_d$. To each of these $\rho$-variables of degree $d$ there
correspond $d$ variables called the {\em roots} of $x_{di}$ or simply
$\rho$-roots and are written as:
$\xi_{di},\xi^q_{di},\cdots,\xi^{d-1}_{di}, i=1,2,\cdots,r_d ; d=1,2\cdots$.
These $\rho$-roots can be thought of as the characteristic roots of a typical
class $c=(f_{11},f_{12},\cdots,f_{1r_1},f_{21},\cdots,f_{2r_2}\cdots,f_{dr_d}
,\cdots)$ of $GL(n,q)$. Let $F$ denote the set of all irreducible monic
polynomials over GF($q$). A mapping  \mbox{$\alpha:X^{\rho}\rightarrow F$}
such that $deg(x \alpha)$ divides $deg(x)$ for each $x \in X^{\rho}$ is
called
a {\em substitution} of $X^{\rho}$ into $F$. We can apply a substitution
$\alpha$ to $\rho$-roots as well. For each $\xi_{di} \in X^{\rho}$
we choose any root $\gamma_{di}$ of $x_{di}\alpha$ and define
\mbox{$xi^{q^u}_{di}\alpha=\gamma^{q^u}_{di}$}.
Two substitutions $\alpha$ and $\alpha^{'}$ of $X^{\rho}$ into $F$ are said
to be {\em equivalent} if there is degree preserving permutation $\phi$ of
$X^{\rho}$ such that $\alpha\circ\phi=\phi^{'}$. A class of equivalent
substitutions is called {\em mode of substitution}.
Let $\alpha$ be a substitution of $X^{\rho}$ into $F$ and let $f\in F$ be
of degree $d$. For each positive integer $i$ let
\mbox{$
        r_{i}(\alpha,f) = \# \{x\in X^{\rho}: deg(x)=id\ and\ x\alpha=f\}.
     $}.
Then we define
\mbox{$\rho(\alpha,f)=1^{r_{1}(\alpha,f)}.2^{r_{2}(\alpha,f)}\cdots$}
It is easily seen that two substitutions $\alpha$ and  $\alpha^{'}$ are
equivalent if and only if $\rho(\alpha,f)=\rho(\alpha^{'},f)$ for all
$f \in F$. So if $m$ denotes the equivalence class of substitutions to
which $\alpha$ belongs, we can write without ambiguity
$\rho(\alpha,f)=\rho(m,f)$. We shall also talk of substitutions of
$X^{\rho}$
into a class $c$ of $GL(5,q)$. If $c= (\cdots f^{\nu_{c}(f)} \cdots)$ is a
conjugate class of $GL(5,q)$ then a substitution $\alpha$ of $X^{\rho}$
into $F$ satisfying $\mid\rho(\alpha,f)\mid=\mid\nu_{c}(f)\mid$ for all
$f\in F$ will be called a {\em substitution of $X^{\rho}$ into the
class $c$}.

\underline{Example}: Let $\rho=1^{3}.2$ and $c=f^{1^2}_{11}f^{1^3}_{12}$
where both $f_{11}$ and $f_{12}$ are of degree one. So in this case we have
$X^{\rho}=\{x_{11},x_{12},x_{13},x_{21}\}$ and there are two modes of
substitutions of $X^{\rho}$ into $c$, the representative substitutions
being:
\[
   \alpha:\left\{
   \begin{array}{lllll}
      x_{11},x_{12}\rightarrow f_{11}  & & & &
                \hspace{0.5in}          \rho(\alpha,f_{11})=1^2\\
      x_{13},x_{21}\rightarrow f_{12}  & & & &
                \hspace{0.5in}           \rho(\alpha,f_{12})=1.2
   \end{array}
   \right.
\]
\[
   \beta:\left\{
   \begin{array}{lllll}
      x_{21}\rightarrow f_{11}  & & & &
                 \hspace{0.2in}     \rho(\beta,f_{11})=2\\
      x_{11},x_{12},x_{13}\rightarrow f_{12}  & & & &
                 \hspace{0.2in}     \rho(\beta,f_{12})=1^3
   \end{array}
   \right.
\]
Substitutions  can also be applied to $\rho$-roots. For each $x_{di} \in
X^{\rho}$ choose any root $\gamma_{di}$ of $x_{di}\alpha$ and define
$\xi^{q^u}_{di}\alpha=\gamma^{q^u}_{di} (u=0,1,2\cdots,d-1)$. For our
purposes it would be irrelevant what root $\gamma_{di}$ is chosen.

\subsection{$\rho$-Functions}

A $\rho$-{\em function}
$U_{\rho}(x^{\rho})=U_{\rho}(x_{11},x_{12},\cdots,x_{1r_1},x_{21},\cdots,
x_{2r_2},\cdots)$
is a function of $\rho$-variables. It is a complex-valued function on the
set of modes of substitutions of $X^{\rho}$ into $F$. It can also be
considered as a function of $\rho$-roots if we write
\[ U_{\rho}(\xi^{\rho})=
   U(\xi_{11},\xi_{12},\cdots,\xi_{1r_1},\xi_{21},\cdots,\xi_{2r_2}.\cdots).
\]
Let $\gamma\in GF(q^d)^*$ be a root of $x_{di}\alpha$. Then
\[ U_{\rho}(\xi^{\rho}m)= U_{\rho}(\xi^{\rho}\alpha)
U(\gamma_{11},\gamma_{12},\cdots,\gamma_{1r_1},\gamma_{21},
\cdots,\gamma_{2r_2}.\cdots).
\]
where $m$ is the mode of the substitution $\alpha$.

\subsection{Uniform Functions}
For each partition $\rho$ of $n$ let there be a given a $\rho$-function
$U_{\rho}(x^{\rho})$. Then a {\em uniform function U on} $GL(n,q)$ is the
class function whose value on the class $c$ is given by:
\[
U(c)=\sum_{\rho}\sum_{m}Q(m,c)U_{\rho}(x^{\rho}m).
\]
Here summation is taken over all partitions $\rho$ of $n$ and all modes
$m$ of substitutions of $X^{\rho}$ into the class $c$ , and
\[
Q(m,c) = \prod_{f\in F}\frac{1}{z_{\rho(m,f)}}
Q_{\rho(m,f)}^{\nu_c(f)}(q^{deg(f)})
\]

The functions $U_{\rho}(x^{\rho})$ are called the {\em principal parts}
of $U$ and $U_{\rho}=U_{\rho}(x^{\rho})$ is its $\rho$-part. A uniform
function $U$ whose principal parts are all zero except for $U_{\rho}$ is
called {\em Basic Uniform Function of type} $\rho$.
Consider the subgroup of $GL(m+n,q)$ consisting of the matrices of the
form:
\[
      \left[ \begin{array}{cc}
                   A_{1} & 0 \\
                   \ast     & A_{2}
             \end{array}
      \right]
\]
where $A_1 \in GL(n,q)$ and $A_2 \in GL(m,q)$. Green [1] has shown that if
$U$ and $V$ are uniform functions on $GL(m,q)$ and $GL(n,q)$ respectively
then the function induced on $GL(m+n,q)$ by the function
$\Psi(A)= U(A_1)V(A_2)$ is a uniform function on $GL(m+n,q)$. This induced
uniform function on $GL(m+n,q)$ will be denoted by $U\circ V$.
\section{Calculations of $Q(m,c)$ And $\chi(m,e)$}
In this section we shall describe the computations of the values of the
functions  $Q(m,c)$ and
$\chi(m,e)$. These functions are in  a way dual to each other with respect
to characters of $GL(n,q)$. Their definition is similar in as much as the
Green Polynomials occurring in the definitions $Q(m,c)$ are replaced by
the characters of the symmetric groups in the definition of $\chi(m,e)$,
whereas, of course, the class $c$ is replaced by the dual class $e$. For
a given partition $\rho$ and for a given class $c$ and dual class $e$
we have:
\[
Q(m,c) = \prod_{f\in F}\frac{1}{z_{\rho(m,f)}}
Q_{\rho(m,f)}^{\nu_c(f)}(q^{deg(f)})
\]
and
\[
\chi(m,e) = \prod_{g\in G}\frac{1}{z_{\rho(m,e)}}
Q_{\rho(m,e)}^{\nu_e(g)}(q^{deg(g)})
\]
As has been pointed out earlier here $m$ is the mode of substitution
of $X^{\rho}$ into the class $c$ and $z_{\rho(m,f)}$ is the order of
centralizer (in the symmetric group) of an element represented by the
partition $\rho(m,f)$. If there are more than one modes of substitutions
we get more than one such functions (as happens in the case of
$\rho = 1^3.2$ and the classes $C, E$).

Given the partition $\rho$ and the class $c$ we first find the mode(s) of
substitutions of $X^{\rho}$ into $c$. For each mode we evaluate the
coorresponding $\rho(m,f)$ and $z_{\rho(m,f)}$. We read the value of
$Q_{\rho(m,f)}^{\nu_c(f)}(q)$ from the table in Green [1] for each
irreducible polynomial appearing in class $c$, and can thus compute the
value of $Q(m,c)$.

\mbox {\bf Examples}
\begin{enumerate}
   \item  Let $\rho=1^5$  and $c=B_1=f_{1a}^{1^4}.f_{1b}^1$ , so
      $X^{\rho} = \{x_{11},x_{12},x_{13},x_{14},x_{15} \}$. There is only
      one mode of substitution of $X^{\rho}$ into $c$ whose representative
      substitution being
      $x_{11},x_{12},x_{13},x_{14} \rightarrow f_{1a}$ and
      $x_{15} \rightarrow f_{1b}$. So
      $\rho(m,f_{1a})=1,
      \rho(m,f_{1a})=1,
      \rho(m,f_{1b})=1,
      z_{\rho(m,f_{1a})}=4.1,
      z_{\rho(m,f_{1b})}=1$. Now

$
      Q(m,c)=
      \frac{1}{z_{\rho(m,f_{1a})}}Q^{\nu_{c}(f_{1a})}_{\rho(m,f_{1a})}(q)
      \frac{1}{z_{\rho(m,f_{1b})}}Q^{\nu_{c}(f_{1b})}_{\rho(m,f_{1b})}(q)
      = \frac{1}{4!}Q^{1^4}_{1^4}(q) Q^{1}_{1}(q)
$

$
            =  \frac{1}{24}(q+1)(q^2+q+1)(q^3+q^2+q+1)
$

   \item  Let $\rho = 2.1^3, c=f^{2.1}_{1a} f^{1^2}_{1b}$. So
       $X^{\rho}= \{ x_{21}, x_{11},x_{12},x_{13}\} $. There are two
       modes of substitutions of  $X^{\rho}$ into $c$ namely:
       The mode $m_1 : x_{11}, x_{12}, x_{13} \rightarrow f_{1a} , \
                 x_{21}\rightarrow f_{1b}$ and the
       mode $m_2:
       x_{11},x_{21} \rightarrow f_{1a},\ x_{12},x_{13} \rightarrow f_{1b}$.
       We have $\rho(m_1,f_{1a})= 1^3,
                \rho(m_1,f_{1b})= 2,
                \rho(m_2,f_{1a})= 1.2,
                \rho(m_2,f_{1b})= 1^2$.
       So $Q(m_1,c)= \frac{1}{3!}\frac{1}{2}Q^{1^3}_{1^3}(q)Q^{1^2}_{2}(q)=
                     - \frac{1}{12}(q+1)(q^3+1)$ and
          $Q(m_2,c)= \frac{1}{2}\frac{1}{2}Q^{1^3}_{1.2}(q)Q^{1^2}_{1^2}(q)=
                     - \frac{1}{4}(q+1)(q^3-1)$.
\end{enumerate}

The determination of the functions $\chi(m,e)$ is similar excepting for the
fact that the class $c$ is replaced by the dual class $e$ and the irreducible
polynomials are replaced by the simplexes $g$ and the Green Polynomials
$Q^{\nu}_{\rho}(q)$ are replaced by the characters $\chi^{\nu}_{\rho}$ of the
symmetric groups. As a matter of fact in most cases the value of
$\chi(m,e)$ is equal to the coefficient of the highest power of $q$
occurring
in the corresponding $Q(m,c)$. The table 3 and table 4 (resp.) give the
values of $Q(m,c)$
and $\chi(m,e)$ (resp.) for all classes $c$ and dual classes $e$, and for
all partitions $\rho_{i}, 1 \leq i \leq 7$ of 5.


\section{Basic characters}
The {\em Basic Characters} are the basic constituents of the irreducible
characters of $GL(5,q)$. The 42 types of irreducible characters of $GL(5,q)$
can be grouped into 17 categories such that the characters in the same
category involve the same basic constituents namely the basic characters.
Roughly speaking these basic characters express the irrationalities
involved in the irreducible characters.

Let $\rho = 1^{r_1}.2^{r_2}\cdots$ be a partition of $n$ and let $h$
denote the vector
$(h_{11},h_{12},\cdots,h_{1r_1},h_{21},\cdots,h_{2r_2},\cdots)$ whose
components (which are integers) correspond to parts of $\rho$. Then
$B^{\rho}(h)$ {\em the basic character of type} $\rho$ is defined to be
the character
\[
J_1(h_{11})\circ J_1(h_{12})\circ \cdots J_1(h_{1r_1})\circ
J_2(h_{21})\circ \cdots J_2(h_{2r_2})\circ \cdots
\]
of $GL(n,q)$ where $J_d(k)$  is a character of $GL(d,q)$ defined as:
$J_d(k)(c) = 0$ if the class $c$ is not primary and
$J_d(k)(f^{\lambda}) = k(\lambda,q^{d(f)}
\{\theta^{k}(\gamma)+  \theta^{kq}(\gamma)+ \theta^{kq^2}(\gamma)+\cdots
                               + \theta^{q^{k(d(f)-1)}}(\gamma) \}
$
where $\gamma$ is a root of the polynomial $f$ and $\theta$ is a
generating character of the the group $GF(q^{n!})^{*}$ and we have
abbreviated $deg(f)$ as $d(f)$. $B^{\rho}(h)$ is a basic uniform
function whose $\rho$-part $B_{\rho}(h:\xi^{\rho})$ is given by the
expression
\[
  \prod_{d} \{\sum S_{d}(h_{d1}:\xi_{d1})\cdots S_{d}(h_{dr_d}:\xi_{dr_d}\}
\]
where the summation is over all permutations of $1,2,\cdots,r_d$ and
\[
  S_{d}(k:\xi) =
   \theta^{k}(\xi)+\theta^{qk}(\xi)+\cdots+\theta^{q^{d-1}k}(\xi).
\]
Actually $\theta^{k}(\xi)$ makes sense only if $\xi$ is a root of an
irreducible polynomial over $GF(q)$. However when a substitution is
applied to $\xi$, a $\rho$-root is transformed into a root of an
irreducible polynomial over $GF(q)$. The value of $B^{\rho}$ at a class
$c$ is given by
\[
   B^{\rho}(h)(c)= \sum_{m} Q(m,c) B_{\rho}(h:\xi^{\rho}m)
\]
where the summation is taken over all modes of substiitutions of $X^{\rho}$
into the class $c$.

\subsection{Simplex and Dual Classes}
Let $f \in F$ be an irreducible polynomial of degree $s$ over $GF(q)$. Let
$\omega$ be a generator of $GF(q^{s})^{*}$. Then a root of $f$ has the form
$\omega^{k}$ where $k$ is uniquely determined modulo($q^{s}-1$), and the set
\{$\omega^{k},\omega^{kq},\cdots,\omega^{kq^{s-1}} $\}
consists all the roots of $f$. Recall that the set
\{$k,kq,kq^2,\cdots,kq^{s-1}$ \} is an $s$-simplex with
$k,kq,kq^2,\cdots,kq^{s-1}$ as its roots and $s$ its degree. It is clear that
there are as many simplexes of degree $s$ as there are irreducible polynomials
$f \in F$ of degree $s$ over $GF(q)$. As before the set of all s-simplexes
for $1\leq s \leq n$ will be denoted by $G$ and a particular simplex by the
letter $g$. For each partition $\rho = 1^{r_1}.2^{r_2}\cdots $, let
$Y^{\rho}$  denote the set of $\rho$-variables
\{$y_{11},y_{12},\cdots,y_{1r_1},y_{21},\cdots,y_{2r_2},\cdots$\}. The
$\rho$-variable $y_{ij}$ is considered to be of degree $i$.
A mapping $\alpha : Y^{\rho} \rightarrow G$ is called a
{\em substitution} of $Y^{\rho}$ into $G$ if $deg(y\alpha)$ divides $deg(y)$.
The notions of equivalenc of two substitutions, mode of a substitution and
the partitions $\rho(\alpha,g)$ are defined exactly in the same manner as
for the substitutions of $X^{\rho}$ into $F$. A function
\mbox {$\nu \rightarrow P$} such that
\mbox {$\sum_{g \in G} \mid\nu(g)\mid deg(g) = n$}
specifies what is called a {\em dual class} and will be denoted by the
symbol $e=(\cdots,g^{\nu(g)},\cdots)$. To each symbol of dual class there
will correspond a character which will be denoted by the same symbol.
For a fixed partition $\rho$ of $n$ let $m$ be the mode of certain
substitution $\alpha$ of $Y^{\rho}$ into the set $G$ of all simplexes. Set
$h_{di}m = c_{di}(q^d-1)/(q^{s_{di}}-1)$, where $c_{di}$ is a root
of the simplex $y_{di}\alpha$ and $s_{di}$ is the degree of $y_{di}\alpha$

We shall now describe the computation of $\rho$-parts  of the
basic characters of $GL(5,q)$ for a partition $\rho$ of $5$ corresponding
to a given symbol $e$ of a dual class of $GL(5,q)$. First we find the modes
of substitution of $Y^{\rho}$ into $e$, then for each mode $m$ we find the
integers $h_{di}m$ defined above. Then the $\rho$-parts of the basic
character for the mode $m$ corresponding to the irreducible character given
by the symbol $e$ of the dual class is given by
\[
   B^{\rho}(h^{\rho}m:\xi^{\rho})=
   B^{\rho}(h_{11}m,h_{12}m,\cdots,h_{1r_1}m, h_{21}m,\cdots,h_{2r_2}m,\cdots
      :\xi^{\rho} )
\]
where now $h_{11}m, h_{12}m, \cdots$ are all integers and so the value of the
above expression can be calculated.

For example let $\rho = 1.2^2, e=g^{4}_{1i}g_{1j}$. Then
$Y^{\rho} = \{ y_{11}, y_{21}, y_{22} \}$. It is easily seen that there is
only one mode of substitution of $Y^{\rho}$ into $e$ and a substitution
$\alpha : Y^{\rho} \rightarrow e $ belonging to this mode is given by
$y_{21}, y_{22} \rightarrow g_{1i}, y_{11} \rightarrow g_{1j}$.
Also $h^{\rho} = \{ h_{11}, h_{21}, h_{22} \}, h_{11}m=j$, where $j$ is a
root of $g_{1i}$, and $h_{21}m = h_{22}m = i(q+1)$ where $i$ is root of
$g_{1i}$. So $B^{\rho}(h^{\rho}m) = B^{\rho}(j,i(q+1), i(q+1))$ and the
$\rho$-part of the basic character corresponding to the symbol
$e= g^{4}_{1i} g_{1j}$ of the dual class is:
\newline
$
   B^{\rho}(h^{\rho}m:\xi^{\rho}) =
   S_{1}(j,\xi_{11}) [ \sum S_{2}(i(q+1), \xi_{21}) S_{2}(i(q+1),\xi_{22}) ]
$
\newline
$
   = 2 S_{1}(j,\xi_{11}) S_{2}(i(q+1), \xi_{21}) S_{2}(i(q+1),\xi_{22})
$

Let us determine the value of this basic character on class $c$ of type
$A$ say $f_{1a}^5$. We first find the modes of substitutions of $X^{\rho}$
into $f_{1a}^5$. As $X^{\rho} = \{x_{11},x_{21},x_{22}\}$ there is only one
mode of substitution whose representative substitution maps each of the
$\rho$-variable into $f_{1a}$. According to our notation the root of the
polynomial can be expressed as $\omega_{1}^{a}$. So we have
\(
   \xi_{11}m= \xi_{21}m=\xi_{22}m=\omega_{1}^{a}
\)
and
\(
   B^{\rho}(h^{\rho}m:\xi^{\rho})(c) =
   2S_{1}(j,\omega_{1}^{a}) (S_{2}(i(q+1),\omega_{1}^{a})^{2}=
   2 \theta^{j} \omega_{1}^{a}
   (\theta^{i(q+1)} \omega_{1}^{a} + \theta^{iq(q+1)} \omega_{1}^{a})^{2}=
   2\varepsilon_{1}^{ja}
   (\varepsilon_{1}^{i(q+1)a} + \varepsilon_{1}^{iq(q+1)a})^{2}  =
   8 \varepsilon_{1}^{(4i+j)a}
\)
,where we recall $\varepsilon_{1}$ is a primitive $(q-1)$-th root of unity.
The value of this basic charcater is same on every class of type $A$ namely
on each of the classes $A_{i}, 1 \leq i \leq 7$.

Table 6 gives the value of the $\rho$-part of each of the $17$ types of
basic characters on each of the $17$ types of classes corresponding to
each partition $\rho$ of $5$ for the $GL(5,q)$. Whenever there is more
than one mode of substitution we have mentioned it explicitly. If the
value of a basic character corresponding to a partition $\rho$ is missing
then it should be understood that there is no mode of substitution - either
into the symbol $e$ of a dual class or into the class $c$.


\section{Irreducible Characters of $GL(5,q)$}

The character of $GL(5,q)$ corresponding to the dual class
$e = (\ldots,g^{\nu(g)},\ldots)$ of $GL(5,q)$ is given Green [1]
by the expression
\[
(-1)^{5-\sum \mid\nu(g)\mid}\sum_{\rho, m} \chi(m,e) B^{\rho}(h^{\rho}m)
\]
where the summation is over all the $7$ partitions $\rho_{i}$ ,
$1 \leq i \leq 7 $ of $5$ and over all modes of substitutuions of
$Y^{\rho{i}}$ into the dual class $e$. The degree of this character is
given by the expression
\[
   \Psi_{5}(q) \prod_{g \in G} (-1)^{\mid \nu(g)\mid} \{ \nu(g):q^{deg(g)} \}
\]
where $\Psi_{n}(q) = \prod_{i=1}^{i=n} (q^i -1)$ and if  $\lambda$
is a partition
of $n$ with parts $\ell_{1} \geq \ell_{2} \geq \cdots \ell_{p} \geq$ then
\[
   \{\lambda:q \} = q^{\ell_2+2\ell_3+\cdots}
   \prod_{1 \leq r < s \leq p}(1- q^{(\ell_r-\ell_s-r+s)})/
   \prod_{r=1}^{r=p}\Phi_{\ell_r+p-r}(q)
\]
where $\Phi_{n}(q) = (-1)^n \Psi_{n}(q)$. Now to calculate the value of
$ \sum_{m}\chi(m,e)B^{\rho}(h^{\rho}m)$, for a specific partition $\rho$,
at a class $c$ of $GL(5,q)$ we write
$ B^{\rho}(h^{\rho}m)$  as simply
$ B^{\rho}(h)$ and keep in mind that now $h=h^{\rho}m$ is a specific
sequence of integers whose value is uniquely determined by the partition
$\rho$ of and the mode $m$ of substitution of $Y^{\rho}$ into $e$. Recall
that the value of $B^{\rho}$ at a class $c$ of $GL(5,q)$ is given by the
expression $\sum_{m}Q(m,c)B_{\rho}(h:\xi^{\rho}m)$, where the summation
over all modes of substitutions of $X^{\rho}$ into the class $c$ and
$B_{\rho}(h:\xi^{\rho}m)$ is the $\rho$-part of the basic character
corresponding to the dual class $e$ whose value has already been
determined (table [6] ). The values of the functions $Q(m,c)$ have
been determined earlier (table [3] ). We do this for each
partition $\rho_{i}, 1 \leq i \leq 7$ of $5$ and for all modes of
substitutions of $Y^{\rho_{i}}$ into the dual class $e$ and further
for all modes of substitutions of $X^{\rho}$ into the class $c$. It may
appear that it is an involved process but in the case of $GL(5,q)$
the computations are managebale. We shall give two example of explicit
computation of values of characters of $GL(5,q)$.

\underline{Examples}

\begin{enumerate}
\item
   Let us determine the value of the character of type $A_{1}^{'}$ which
corresponds to the dual class $e=g^{1^5}_{1i}, 1 \leq i \leq q-1$ on
the class $c = A_1$ (class symbol $f^{1^5}_{1a}, 1 \leq a \leq q-1)$. For
each $\rho_{i}, 1 \leq i \leq 7$, there is only one mode of substitution
of $Y^{\rho_{i}}$ into the dual class $e$ namely each variable mapping
into $g_{1i}$.
Similarly for each $\rho_{i}, 1 \leq i \leq 7$, there is only one mode
of substitution of $X^{\rho_{i}}$ into the class $c$. So in this case
for each $\rho_{i}, 1 \leq i \leq 7$, we find the value of $\chi(m,e)$
corresponding to the dual class $e=A_{1}^{'}$, multiply it with the value
of $Q(m,c)$ for corresponding $\rho_{i}$ and $c=A_1$, and with the
$\rho_{i}$-part of the basic character of type $A^{'}$ and add these up.
Thus we get the following value for this character on given class:
\[
\frac{1}{120}(q+1)(q^2+q+1)(q^3+q^2+q+1)(q^4+q^3+q^2+q+1)\varepsilon_1^{5ia}+
\frac{1}{12}(q^2+q+1)(q^3+q^2+q+1)(q^5-1)\varepsilon_1^{5ia}+
\]
\[
\frac{1}{8}(q^2+1)(q^3-1)(q^5-1)\varepsilon_1^{5ia}+
\frac{1}{6}(q+1)(q^4-1)(q^5-1)\varepsilon_1^{5ia}+
\frac{1}{6}(q-1)(q^4-1)(q^5-1)\varepsilon_1^{5ia}+
\]

$
\frac{1}{4}(q^2-1)(q^3-1)(q^5-1)\varepsilon_1^{5ia}+
\frac{1}{5}(q-1)(q^2-1)(q^3-1)(q^4-1)\varepsilon_1^{5ia}
$

which when simplified reduces to $q^{10}\varepsilon_{1}^{51a}$. The degree
of each of these characters (as computed by the degree formula mentioned
earlier) is $q^{10}$
\item
   Let us determine the value of the character of type $C_{1}^{'}$ on the
   classes of type $E_2$. In this case there are two modes of substitutions
   for $\rho=\rho_{2}= 2.1^3$ into the dual class of type $C_1^{'}$ as well
   as into the class of type $E^{'}$. Also corresponding to $\rho=\rho_6,
   \rho_7$ there are no modes of substitutions of $Y^{\rho}$ into $C_1^{1}$,
   and corresponding to $\rho=\rho_4, \rho_5, \rho_6, \rho_7$ there are no
   modes of substitutions of $X^{\rho}$ into the class $E_2$ .
   Thus corresponding to the partitions $\rho_1$ and $\rho_3$ the
   contribution to the value of character is:
\[
   (q+1) [ \frac{1}{2}\varepsilon_{1}^{i(2a+b)+j(a+b)}+
           \frac{1}{4}\varepsilon_{1}^{i(2a+c)+2jb}   +
           \frac{1}{2}\varepsilon_{1}^{i(a+2b)+j(a+c)}   +
           \varepsilon_{1}^{i(a+b+c)+j(a+b)}   +
           \frac{1}{4}\varepsilon_{1}^{i(2b+c)+2ja} ]
\]

$
           -\frac{1}{4}(q-1)[\varepsilon_{1}^{i(2a+c)+2jb} +
            \varepsilon_{1}^{i(2b+c)+2ja} ]
$

To find the contribution corresponding to the partitions $\rho_2$ we have
to take into consideration both the modes of substitutions. This
contributions comes out to be
\newline
$
   -\frac{1}{4}(q+1)\varepsilon_{1}^{i(2a+c)+2jb} +
    \frac{1}{4}(q-1)\varepsilon_{1}^{i(2b+c)+2ja}
   -\frac{1}{4}(q+1) [2\varepsilon_{1}^{i(a+2b)+j(a+c)} +
                        \varepsilon_{1}^{i(2b+c)+2ja} ]
$
\newline
$
   +\frac{1}{4}(q-1) [2\varepsilon_{1}^{i(2a+b)+j(b+c)} +
                        \varepsilon_{1}^{i(2a+c)+2jb} ]
$

Adding the two contributions we get the value of the characters of type
$C_1^{'}$ on the classse of type $E_2$ (as can be seen there is no
contribution from the partitions $\rho_4, \rho_5, \rho_6, \rho_7$).
This sum comes out to be equal to
\[
   q \varepsilon_{1}^{i(2a+b)+j(b+c)} +
   (q+1)\varepsilon_{1}^{i(a+b+c)+j(a+b)}
\]
The degree of each of these characters as computed by the degree formula is
$q^4(q^2+1)(q^4+q^3+q^2+q+1)$
\end{enumerate}
A complete listing of characters of $GL(5,q)$ is given in table 7.
%
%












%
\newpage
\begin{center}
   {\bf Appendix} \\
   {\bf List of Tables}
\end{center}
\begin{tabular}{ll}
   Table 1  &  Canonical forms of representatives of conjugate    \\
            &  classes of $GL(5,q)$                 \\
   Table 2  &  Description of Conjugate Classes   \\
   Table 3  &  The functions $Q(m,c)$.             \\
   Table 4  &  The functions $\chi(m,e)$            \\
   Table 5  &  Description of Dual Classes           \\
            &  (Types of Characters)                 \\
   Table 6  &  Basic Characters                      \\
   Table 7  &  Irreducble Characters of $GL(5,q)$     \\
\end{tabular}

\newpage
\begin{center}
   {\bf TABLE \# 1} \\
   {\bf REPRESENTATIVES OF CONJUGATE CLASSES OF GL(5,q)} \\
\end{center}

\begin{displaymath}
A_{1} = \left[
   \begin{array}{ccccc}
              \omega_{1}^{a} & & & & \\
               & \omega_{1}^{a} & & & \\
               & & \omega_{1}^{a} &  & \\
               & & & \omega_{1}^{a}  & \\
               & & & & \omega_{1}^{a}   \\
   \end{array}
        \right]
A_{2} = \left[
   \begin{array}{ccccc}
              \omega_{1}^{a} & & & & \\
               1 & \omega_{1}^{a} & & & \\
               & & \omega_{1}^{a} &  & \\
               & & & \omega_{1}^{a}  & \\
               & & & & \omega_{1}^{a}   \\
   \end{array}
        \right]
A_{3} = \left[
   \begin{array}{ccccc}
              \omega_{1}^{a} & & & & \\
               1 & \omega_{1}^{a} & & & \\
               & & \omega_{1}^{a} &  & \\
               & & 1 & \omega_{1}^{a}  & \\
               & & & & \omega_{1}^{a}   \\
   \end{array}
        \right]
\end{displaymath}
\begin{displaymath}
A_{4} = \left[
   \begin{array}{ccccc}
              \omega_{1}^{a} & & & & \\
               1 & \omega_{1}^{a} & & & \\
               & 1 & \omega_{1}^{a} &  & \\
               & & & \omega_{1}^{a}  & \\
               & & & & \omega_{1}^{a}  \\
   \end{array}
        \right]
A_{5} = \left[
   \begin{array}{ccccc}
              \omega_{1}^{a} & & & & \\
               1 &  \omega_{1}^{a} & & & \\
               & 1 & \omega_{1}^{a} &  & \\
               & & & \omega_{1}^{a}  & \\
               & & & 1 & \omega_{1}^{a}   \\
   \end{array}
        \right]
A_{6} = \left[
   \begin{array}{ccccc}
              \omega_{1}^{a} & & & &  \\
               1 & \omega_{1}^{a} & & & \\
               & 1 & \omega_{1}^{a} &  & \\
               & &  1  & \omega_{1}^{a} & \\
               & & & & \omega_{1}^{a}   \\
   \end{array}
        \right]
\end{displaymath}
\begin{displaymath}
A_{7} = \left[
   \begin{array}{ccccc}
              \omega_{1}^{a} & & & & \\
               1 & \omega_{1}^{a} & & & \\
               & 1 & \omega_{1}^{a} &  & \\
               & & 1 & \omega_{1}^{a}  & \\
               & & & 1 & \omega_{1}^{a}  \\
   \end{array}
        \right]
B_{1} = \left[
   \begin{array}{ccccc}
              \omega_{1}^{a} & & & & \\
               & \omega_{1}^{a} & & & \\
               & & \omega_{1}^{a} &  & \\
               & & & \omega_{1}^{a}  & \\
               & & & & \omega_{1}^{b}   \\
   \end{array}
        \right]
B_{2} = \left[
   \begin{array}{ccccc}
              \omega_{1}^{a} & & & & \\
               1 & \omega_{1}^{a} & & & \\
               & & \omega_{1}^{a} &  & \\
               & & & \omega_{1}^{a}  & \\
               & & & & \omega_{1}^{b}   \\
   \end{array}
        \right]
\end{displaymath}
\begin{displaymath}
B_{3} = \left[
   \begin{array}{ccccc}
              \omega_{1}^{a} & & & & \\
               1 & \omega_{1}^{a} & & & \\
               & & \omega_{1}^{a} &  & \\
               & & 1 & \omega_{1}^{a}  & \\
               & & & & \omega_{1}^{b}   \\
   \end{array}
        \right]
B_{4} = \left[
   \begin{array}{ccccc}
              \omega_{1}^{a} & & & & \\
               1 & \omega_{1}^{a} & & & \\
               & 1 & \omega_{1}^{a} &  & \\
               & & & \omega_{1}^{a}  &  \\
               & & & & \omega_{1}^{b}    \\
   \end{array}
        \right]
B_{5} = \left[
   \begin{array}{ccccc}
              \omega_{1}^{a} & & & & \\
               1 & \omega_{1}^{a} & & & \\
               & 1 & \omega_{1}^{a} &  & \\
               & & 1 & \omega_{1}^{a}  & \\
               & & & & \omega_{1}^{b}   \\
   \end{array}
        \right]
 \end{displaymath}
\begin{displaymath}
C_{1} = \left[
   \begin{array}{ccccc}
              \omega_{1}^{a} & & & & \\
               & \omega_{1}^{a} & & & \\
               & & \omega_{1}^{a} &  & \\
               & & & \omega_{1}^{b}  & \\
               & & & & \omega_{1}^{b}   \\
   \end{array}
        \right]
C_{2} = \left[
   \begin{array}{ccccc}
              \omega_{1}^{a} & & & & \\
               1 & \omega_{1}^{a} & & & \\
               & & \omega_{1}^{a} &  & \\
               & & & \omega_{1}^{b}  & \\
               & & & & \omega_{1}^{b}   \\
   \end{array}
        \right]
C_{3} = \left[
   \begin{array}{ccccc}
              \omega_{1}^{a} & & & & \\
               1 & \omega_{1}^{a} & & & \\
               & 1 & \omega_{1}^{a} &  & \\
               & & & \omega_{1}^{b}  & \\
               & & & & \omega_{1}^{b}   \\
   \end{array}
        \right]
\end{displaymath}
\begin{displaymath}
C_{4} = \left[
   \begin{array}{ccccc}
              \omega_{1}^{a} & & & & \\
               & \omega_{1}^{a} & & & \\
               & & \omega_{1}^{a} &  & \\
               & & & \omega_{1}^{b}  & \\
               & & & 1 & \omega_{1}^{b}   \\
   \end{array}
        \right]
C_{5} = \left[
   \begin{array}{ccccc}
              \omega_{1}^{a} & & & & \\
               1 & \omega_{1}^{a} & & & \\
               & & \omega_{1}^{a} &  & \\
               & & & \omega_{1}^{b}  & \\
               & & & 1 & \omega_{1}^{b}   \\
   \end{array}
        \right]
C_{6} = \left[
   \begin{array}{ccccc}
              \omega_{1}^{a} & & & & \\
               1 & \omega_{1}^{a} & & & \\
               & 1 & \omega_{1}^{a} &  & \\
               & & & \omega_{1}^{b}  & \\
               & & & 1 & \omega_{1}^{b}   \\
   \end{array}
        \right]
\end{displaymath}
%
%
\begin{displaymath}
D_{1} = \left[
   \begin{array}{ccccc}
              \omega_{1}^{a} & & & & \\
               & \omega_{1}^{a} & & & \\
               & & \omega_{1}^{a} &  & \\
               & & & \omega_{1}^{b}  & \\
               & & & & \omega_{1}^{c}   \\
   \end{array}
        \right]
D_{2} = \left[
   \begin{array}{ccccc}
              \omega_{1}^{a} & & & & \\
               1 & \omega_{1}^{a} & & & \\
               & & \omega_{1}^{a} &  & \\
               & & & \omega_{1}^{b}  & \\
               & & & & \omega_{1}^{c}   \\
   \end{array}
        \right]
D_{3} = \left[
   \begin{array}{ccccc}
              \omega_{1}^{a} & & & & \\
               1 & \omega_{1}^{a} & & & \\
               & 1 & \omega_{1}^{a} &  & \\
               & & & \omega_{1}^{b}  & \\
               & & & & \omega_{1}^{c}   \\
   \end{array}
        \right]
\end{displaymath}
%
%
%
\begin{displaymath}
E_{1} = \left[
   \begin{array}{ccccc}
              \omega_{1}^{a} & & & & \\
               & \omega_{1}^{a} & & & \\
               & & \omega_{1}^{b} &  & \\
               & & & \omega_{1}^{b}  & \\
               & & & & \omega_{1}^{c}   \\
   \end{array}
        \right]
E_{2} = \left[
   \begin{array}{ccccc}
              \omega_{1}^{a} & & & & \\
               1 & \omega_{1}^{a} & & & \\
               & & \omega_{1}^{b} &  & \\
               & & & \omega_{1}^{b}  & \\
               & & & & \omega_{1}^{c}   \\
   \end{array}
        \right]
E_{3} = \left[
   \begin{array}{ccccc}
              \omega_{1}^{a} & & & & \\
               1 & \omega_{1}^{a} & & & \\
               &  & \omega_{1}^{b} &  & \\
               & & 1 & \omega_{1}^{b}  & \\
               & & & & \omega_{1}^{c}   \\
   \end{array}
        \right]
\end{displaymath}
%
%
%
\begin{displaymath}
F_{1} = \left[
   \begin{array}{ccccc}
              \omega_{1}^{a} & & & & \\
               & \omega_{1}^{a} & & & \\
               & & \omega_{1}^{b} &  & \\
               & & & \omega_{1}^{c}  & \\
               & & & & \omega_{1}^{d}   \\
   \end{array}
        \right]
F_{2} = \left[
   \begin{array}{ccccc}
              \omega_{1}^{a} & & & & \\
               1 & \omega_{1}^{a} & & & \\
               & & \omega_{1}^{b} &  & \\
               & & & \omega_{1}^{c}  & \\
               & & & & \omega_{1}^{d}   \\
   \end{array}
        \right]
G = \left[
   \begin{array}{ccccc}
              \omega_{1}^{a} & & & & \\
               & \omega_{1}^{b} & & & \\
               &  & \omega_{1}^{c} &  & \\
               & &  & \omega_{1}^{d}  & \\
               & & & & \omega_{1}^{e}   \\
   \end{array}
        \right]
\end{displaymath}
%
%
%
\begin{displaymath}
H_{1} = \left[
   \begin{array}{ccccc}
              \omega_{1}^{a} & & & & \\
               & \omega_{1}^{a} & & & \\
               & & \omega_{1}^{a} &  & \\
               & & & \omega_{2}^{b}  & \\
               & & & & \omega_{2}^{bq}   \\
   \end{array}
        \right]
H_{2} = \left[
   \begin{array}{ccccc}
              \omega_{1}^{a} & & & & \\
               1 & \omega_{1}^{a} & & & \\
               & & \omega_{1}^{a} &  & \\
               & & & \omega_{2}^{b}  & \\
               & & & & \omega_{2}^{bq}   \\
   \end{array}
        \right]
H_{3} = \left[
   \begin{array}{ccccc}
              \omega_{1}^{a} & & & & \\
               1 & \omega_{1}^{a} & & & \\
               & 1 & \omega_{1}^{a} &  & \\
               & & & \omega_{2}^{b}  & \\
               & & & & \omega_{2}^{bq}   \\
   \end{array}
        \right]
\end{displaymath}
%
%
%
\begin{displaymath}
I_{1} = \left[
   \begin{array}{ccccc}
              \omega_{1}^{a} & & & & \\
               & \omega_{1}^{a} & & & \\
               & & \omega_{1}^{b} &  & \\
               & & & \omega_{2}^{c}  & \\
               & & & & \omega_{2}^{cq}   \\
   \end{array}
        \right]
I_{2} = \left[
   \begin{array}{ccccc}
              \omega_{1}^{a} & & & & \\
               1 & \omega_{1}^{a} & & & \\
               & & \omega_{1}^{b} &  & \\
               & & & \omega_{2}^{c}  & \\
               & & & & \omega_{2}^{cq}   \\
   \end{array}
        \right]
J = \left[
   \begin{array}{ccccc}
              \omega_{1}^{a} & & & & \\
               & \omega_{1}^{b} & & & \\
               &  & \omega_{1}^{c} &  & \\
               & & & \omega_{2}^{d}  & \\
               & & & & \omega_{2}^{dq}   \\
   \end{array}
        \right]
\end{displaymath}
%
%
%
\begin{displaymath}
K_{1} = \left[
   \begin{array}{ccccc}
              \omega_{1}^{a} & & & & \\
               & \omega_{2}^{b} & & & \\
               & & \omega_{2}^{bq} &  & \\
               & & & \omega_{2}^{b}  & \\
               & & & & \omega_{2}^{bq}   \\
   \end{array}
        \right]
K_{2} = \left[
   \begin{array}{ccccc}
              \omega_{1}^{a} & & & & \\
               1 & \omega_{2}^{b} & & & \\
               & & \omega_{2}^{bq} &  & \\
               & & 1 & \omega_{2}^{b}  & \\
               & & & 1 & \omega_{2}^{bq}   \\
   \end{array}
        \right]
L = \left[
   \begin{array}{ccccc}
              \omega_{1}^{a} & & & & \\
               & \omega_{2}^{b} & & & \\
               &  & \omega_{2}^{bq} &  & \\
               & & & \omega_{2}^{c}  & \\
               & & & & \omega_{2}^{cq}   \\
   \end{array}
        \right]
\end{displaymath}
%
%
%
\begin{displaymath}
M_{1} = \left[
   \begin{array}{ccccc}
              \omega_{1}^{a} & & & & \\
               & \omega_{1}^{a} & & & \\
               & & \omega_{3}^{b} &  & \\
               & & & \omega_{3}^{bq}  & \\
               & & & & \omega_{3}^{bq^2}   \\
   \end{array}
        \right]
M_{2} = \left[
   \begin{array}{ccccc}
              \omega_{1}^{a} & & & & \\
               1 & \omega_{1}^{a} & & & \\
               & & \omega_{3}^{b} &  & \\
               & & 1 & \omega_{3}^{bq}  & \\
               & & & 1 & \omega_{3}^{bq^2}   \\
   \end{array}
        \right]
N = \left[
   \begin{array}{ccccc}
              \omega_{1}^{a} & & & & \\
               & \omega_{1}^{b} & & & \\
               &  & \omega_{3}^{c} &  & \\
               & & & \omega_{3}^{cq}  & \\
               & & & & \omega_{3}^{cq^2}   \\
   \end{array}
        \right]
\end{displaymath}
%
%
%
\begin{displaymath}
O = \left[
   \begin{array}{ccccc}
              \omega_{2}^{a} & & & & \\
               & \omega_{2}^{aq} & & & \\
               & & \omega_{3}^{b} &  & \\
               & & & \omega_{3}^{bq}  & \\
               & & & & \omega_{3}^{bq^2}   \\
   \end{array}
        \right]
P  = \left[
   \begin{array}{ccccc}
              \omega_{1}^{a} & & & & \\
               & \omega_{4}^{b} & & & \\
               & & \omega_{4}^{bq} &  & \\
               & &  & \omega_{4}^{bq^2}  & \\
               & & &  & \omega_{4}^{bq^3}   \\
   \end{array}
        \right]
Q = \left[
   \begin{array}{ccccc}
              \omega_{5}^{a} & & & & \\
               & \omega_{5}^{aq} & & & \\
               &  & \omega_{5}^{aq^2} &  & \\
               & & & \omega_{5}^{aq^3}  & \\
               & & & & \omega_{5}^{aq^4}   \\
   \end{array}
        \right]
\end{displaymath}
%
\begin{landscape}
%
\newpage
\begin{center}
   Table 2. CONJUGATE CLASSES.
\end{center}
\begin{center}
\begin{tabular}{||l|l|l|l|l||}  \hline
  Type      & Symbol & Range of Values & Number of & Number of Elements in \\
            &        & of a,b,c,d,e    & Classes   & Each Type            \\
                                                       \hline
 $A_{1}$ & $f_{1a}^{1^5}$ & $1\leq a\leq (q-1)$ & $(q-1)$
          & $1$                                                     \\
$A_{2}$ & $f_{1a}^{2.1^3}$ &  same    & same
          & $(q+1)(q^2+1)(q^5-1)$                                     \\
$A_{3}$ & $f_{1a}^{2^{2}.1}$ & same   & same
          & $q(q^2+q+1)(q^4-1)(q^5-1)$                                  \\
$A_{4}$ & $f_{1a}^{3.1^{2}}$ & same   & same
          & $q^{3}(q+1)(q^3-1)(q^4-1)(q^5-1)$                       \\
$A_{5}$ & $f_{1a}^{3.2}$ & same       & same
          & $q^{3}(q^2+q+1)(q^4-1)(q^5-1)$                         \\
$A_{6}$ & $f_{1a}^{4.1}$ & same       & same
          & $q^{5}(q+1)(q^3-1)(q^4-1)(q^5-1)$                    \\
\hline
\end{tabular}
\end{center}
\newpage
\begin{small}
\begin{center}
   Table 2. CONJUGATE CLASSES (Cont'd).
\end{center}
%
%
\begin{tabular}{||l|l|l|l|l||}  \hline
  Type      & Symbol & Range of Values & Number of & Number of Elements in \\
            &        & of a,b,c,d,e    & Classes   & Each Type        \\
                                         \hline
$D_{1}$ & $f_{1a}^{1^3}.f_{1b}^{1}.f_{1c}^{1}$ & $1\leq a,b,c \leq (q-1)$
          & $\frac{1}{2}(q-1)(q-2)(q-3)$ &
          $q^{7}(q^3+q^2+q^+1)(q^4+q^3+q^2+q+1)$                \\
        &         & $a,b,c$ all &      &                            \\
        &         & distinct and &      &                            \\
        &         & $(a,b,c)=(a,c,b)$ &      &                            \\
$D_{2}$ & $f_{1a}^{2.1}.f_{1b}^{1}.f_{1c}^{1}$ & same
          & same &
          $q^{7}(q+1)^{2}(q^2+1)(q^3+q^2+q^+1)(q^5-1))$                \\
$D_{3}$ & $f_{1a}^{3}.f_{1b}^{1}.f_{1c}^{1}$ & same
          & same &
          $q^{8}(q+1)(q^2+q+1)(q^4-1)(q^5-1)$                \\
\hline
\end{tabular}
\end{small}
\newpage
\begin{center}
   Table 2. CONJUGATE CLASSES (Cont'd).
\end{center}
\begin{center}
\begin{tabular}{||l|l|l|l|l||}  \hline
  Type      & Symbol & Range of Values & Number of & Number of Elements in \\
            &        & of a,b,c,d,e    & Classes   & Each Type        \\
                                           \hline
$H_{1}$ & $f_{1a}^{1^3}.f_{2b}^{1} $ &
              $1\leq a \leq (q-1)$, $1\leq b\leq(q^2-2) $
          & $\frac{1}{2}q(q-1)^2 $ &
          $q^{7}(q^2+1)(q^5-1)$                \\
        &         & $b\not\equiv 0 \bmod (q+1) $ all distinct &     &                       \\
        &         & $b\sim bq \bmod(q^2-1) $  &    &                   \\
$H_{2}$ & $f_{1a}^{2.1}.f_{2b}^{1} $ &
               same
          &  same   &
          $q^{7}(q^2+q+1)(q^4-1)(q^5-1)$                \\
$H_{3}$ & $f_{1a}^{3}.f_{2b}^{1} $ &
               same
          &  same   &
          $q^8(q^3-1)(q^4-1)(q^5-1)$                \\
$I_{1}$ & $f_{1a}^{1^2}.f_{1b}^{1}.f_{2c}^{1} $ &
              $1\leq a,b \leq (q-1), a\neq b $
          & $\frac{1}{2}q(q-1)^{2}(q-2) $ &
          $q^{9}(q^2+1)(q^2+q+1)(q^5-1)$                \\
\hline
\end{tabular}
\end{center}
\newpage
\begin{center}
   Table 2. CONJUGATE CLASSES (cont'd).
\end{center}
\begin{center}
\begin{tabular}{||l|l|l|l|l||}  \hline
  Type      & Symbol & Range of Values & Number of & Number of Elements in \\
            &        & of a,b,c,d,e    & Classes   & Each Type              \\
                                                       \hline
$M_{1}$ & $f_{1a}^{1^2}.f_{3b}^{1} $ &
              $1\leq a \leq (q-1), 1\leq b \leq (q^3-2) $   &
            $\frac{1}{3}q(q-1)^{2}(q+1) $ &
          $q^{9}(q^4-1)(q^5-1)$                \\
        &         & $ b\not\equiv 0 \bmod (q^2+q+1) $ &     &                       \\
        &         & $b\sim bq \sim bq^2 \bmod(q^3-1) $  &    &                   \\
$M_{2}$ & $f_{1a}^{2}.f_{3b}^{1} $ &
                same               &
                same               &
          $ q^{9}(q^2-1)(q^4-1)(q^5-1)$                \\
%
$N$ & $f_{1a}^{1}.f_{1b}^{1}.f_{3c}^{1} $ &
              $1\leq a,b \leq (q-1), a\neq b   $   &
            $\frac{1}{6}q(q-1)^{2}(q+1)(q-2) $ &
          $q^{10}(q+1)(q^4-1)(q^5-1)$                \\
        &         & $(a,b)=(b,a), 1\leq c\leq (q^3-2) $ &   &      \\
\hline
\end{tabular}
\end{center}
\end{landscape}

\begin{landscape}
%
\newcommand{\ra}{\mbox{$\rho_{1}$}}
\newcommand{\rb}{\mbox{$\rho_{2}$}}
\newcommand{\rc}{\mbox{$\rho_{3}$}}
\newcommand{\rd}{\mbox{$\rho_{4}$}}
\newcommand{\re}{\mbox{$\rho_{5}$}}
\newcommand{\rf}{\mbox{$\rho_{6}$}}
\newcommand{\rg}{\mbox{$\rho_{7}$}}
\newpage
\begin{center}
   Table 3. THE FUNCTIONS $Q(m,c)$.
\end{center}
\begin{center}
\begin{tabular}{||l|l|l|l||}  \hline
  Partition     & Class $A_{1}$ & Class $A_{2}$ & Class $A_{3}$   \\
                                       \hline
\ra\  &  $\frac{1}{120}(q+1)(q^3-1)(q^4-1)(q^5-1)/(q-1)^3    $
        &  $\frac{1}{120}(q+1)(q^2+q+1)(4q^3+3q^2+2q+1)        $
        &  $\frac{1}{120}(q+1)(5q^3+6q^2+3q+1)                 $  \\
\rb\  &  $-\frac{1}{12}(q^2+q+1)(q^3+q^2+q+1)(q^5-1)         $
        &  $\frac{1}{12}(q^2+q+1)(-2q^4+q^3+q^2+1)             $
        &  $\frac{1}{12}(q+1)(-q^3+2q^2+q+1)                   $  \\
\rc\  &  $\frac{1}{8}(q^2+1)(q^3-1)(q^5-1)         $
        &  $-\frac{1}{8}(q^2+1)(q^3-1)               $
        &  $\frac{1}{8}(q^4-q^3+q^2+1)               $  \\
\rd\  &  $\frac{1}{6}(q+1)(q^4-1)(q^5-1)         $
        &  $-\frac{1}{6}(q^2-1)(-q^4+q^3+q^2+q+1)               $
        &  $-\frac{1}{6}(q+1)(q^3-1)               $  \\
\re\  &  $-\frac{1}{6}(q-1)(q^4-1)(q^5-1)         $
        &  $\frac{1}{6}(q-1)(q^5-1)                 $
        &  $-\frac{1}{6}(q-1)(q^3+1)               $  \\
\rf\  &  $-\frac{1}{4}(q^2-1)(q^3-1)(q^5-1)         $
        &  $\frac{1}{4}(q^2-1)(q^3-1)                 $
        &  $-\frac{1}{4}(q-1)(-q^3+q+1)               $  \\
\rg\  &  $\frac{1}{5}(q-1)(q^2-1)(q^3-1)(q^4-1)         $
        &  $-\frac{1}{5}(q-1)(q^2-1)(q^3-1)                 $
        &  $\frac{1}{5}(q-1)(q^2-1)               $  \\
                           \hline
\end{tabular}
\end{center}
%
\vspace{0.5in}
\begin{center}
\begin{tabular}{||l|l|l|l|l||}  \hline
  Partition     & Class $A_{4}$ & Class $A_{5}$ & Class $A_{6}$
                & Class $A_{7}$                                    \\
                                       \hline
 $\ra$  &  $\frac{1}{120}(q+1)(6q^2+3q+1)                      $
        &  $\frac{1}{120}(5q^2+4q+1)                           $
        &  $\frac{1}{120}(4q+1)                                $
        &  $\frac{1}{120}                                      $  \\
 $\rb$  &  $\frac{1}{12}(3q^2+q+1)                             $
        &  $\frac{1}{12}(q+1)^2                                $
        &  $\frac{1}{12}(2q+1)                                 $
        &  $\frac{1}{12}                                       $  \\
 $\rc$   &  $\frac{1}{8}(-2q^3+q^2+1)                           $
        &  $\frac{1}{8}(q^2+1)                                 $
        &  $\frac{1}{8}                                        $
        &  $\frac{1}{8}                                        $  \\
 $\rd$   &  $\frac{1}{6}(q+1)                                   $
        &  $\frac{1}{6}(-q^2+q+1)                              $
        &  $\frac{1}{6}(q+1)                                   $
        &  $\frac{1}{6}                                        $  \\
 $\re$   &  $-\frac{1}{6}(q-1)                                  $
        &  $\frac{1}{6}(q^2-q+1)                               $
        &  $\frac{1}{6}(q-1)                                   $
        &  $\frac{1}{6}                                        $  \\
 $\rf$   &  $-\frac{1}{4}(q^2-1)                              $
        &  $-\frac{1}{4}(q^2-1)                              $
        &  $\frac{1}{4}                                      $
        &  $\frac{1}{4}                                      $  \\
\hline
\end{tabular}
\end{center}
\newpage
\begin{center}
   Table 3. THE FUNCTIONS $Q(m,c)$  (cont'd).
\end{center}
\begin{center}
\begin{tabular}{||l|l|l|l|l|l||}  \hline
  Partition     & Class $B_{1}$ & Class $B_{2}$ & Class $B_{3}$
                & Class $B_{4}$ & Class $B_{5}$                     \\
                                       \hline
 $\ra$  &  $\frac{1}{24}(q+1)(q^2+q+1)(q^3+q^2+q+1)           $
        &  $\frac{1}{24}(q+1)(3q^2+2q+1)                      $
        &  $\frac{1}{24}(q+1)(2q+1)                           $
        &  $\frac{1}{24}(3q+1)                                $
        &  $\frac{1}{24}                                       $   \\
 $\rb$  &  $-\frac{1}{4}(q^2+q+1)(q^4-1)                     $
        &  $\frac{1}{4}(-q^3+q^2+q+1)                        $
        &  $\frac{1}{4}(q+1)                                 $
        &  $\frac{1}{4}(q+1)                                 $
        &  $\frac{1}{4}                                      $   \\
 $\rc$  &  $\frac{1}{8}(q-1)(q^2+1)(q^3-1)                     $
        &  $-\frac{1}{8}(q-1)(q^2+1)                          $
        &  $\frac{1}{8}(2q^2-q+1)                                 $
        &  $-\frac{1}{8}(q-1)                                 $
        &  $\frac{1}{8}                                      $   \\
 $\rd$  &  $\frac{1}{3}(q^2-1)(q^4-1)                        $
        &  $-\frac{1}{3}(q^2-1)                          $
        &  $-\frac{1}{3}(q^2-1)                                 $
        &  $\frac{1}{3}                                    $
        &  $\frac{1}{3}                                      $   \\
 $\re$  &  $-\frac{1}{4}(q-1)(q^2-1)(q^3-1)                        $
        &  $\frac{1}{4}(q-1)(q^2-1)                          $
        &  $-\frac{1}{4}(q-1)                                 $
        &  $-\frac{1}{4}(q-1)                                    $
        &  $\frac{1}{4}                                      $   \\
\hline
\end{tabular}
\end{center}

\vspace{0.5in}
\begin{center}
\begin{tabular}{||l|l|l|l|l|l|l||}  \hline
  Partition     & Class $C_{1}$ & Class $C_{2}$ & Class $C_{3}$
                & Class $C_{4}$ & Class $C_{5}$ & Class $C_{6}$        \\
                                       \hline
 $\ra$  &  $\frac{1}{12}(q+1)^2(q^2+q+1)                              $
        &  $\frac{1}{12}(q+1)(2q+1)                                   $
        &  $\frac{1}{12}(q+1)                                         $
        &  $\frac{1}{12}(q+1)(q^2+q+1)                                $
        &  $\frac{1}{12}(2q+1)                                        $
        &  $\frac{1}{12}                                      $   \\
 $\rb$  &  $m_{1}:\ -\frac{1}{12}(q^2-1)(q^2+q+1)                     $
        &  $-\frac{1}{12}(q-1)(2q+1)                                   $
        &  $-\frac{1}{12}(q-1)                                         $
        &  $-\frac{1}{12}(q+1)(q^2+q+1)                                $
        &  $\frac{1}{12}(2q+1)                                        $
        &  $\frac{1}{12}                                              $   \\
        &  $m_{2}:\ -\frac{1}{4}(q+1)(q^3-1)                           $
        &  $-\frac{1}{4}(q+1)                                   $
        &  $-\frac{1}{4}(q+1)                                         $
        &  $-\frac{1}{4}(q^3-1)                                 $
        &  $\frac{1}{4}                                       $
        &  $\frac{1}{4}                                              $   \\
 $\rc$  &  $-\frac{1}{4}(q-1)(q^3-1)                                  $
        &  $-\frac{1}{4}(q-1)                                         $
        &  $-\frac{1}{4}(q-1)                                         $
        &  $-\frac{1}{4}(q+1)(q^3-1)                                $
        &  $\frac{1}{4}                                         $
        &  $\frac{1}{4}                                               $   \\
\hline
\end{tabular}
\end{center}
\newpage
\begin{center}
   Table 3. THE FUNCTIONS $Q(m,c)$  (cont'd).
\end{center}
\begin{center}
\begin{tabular}{||l|l|l|l|l|l|l|l|l|l||}  \hline
  Partition     & Class $D_{1}$ & Class $D_{2}$ & Class $D_{3}$
                & Class $E_{1}$ & Class $E_{2}$ & Class $E_{3}$
                & Class $F_{1}$ & Class $F_{2}$ & Class $G$       \\
                                       \hline
 $\ra$  &  $\frac{1}{6}(q+1)(q^2+q+1)                              $
        &  $\frac{1}{6}(2q+1)                                   $
        &  $\frac{1}{6}                                          $
        &  $\frac{1}{4}(q+1)^2                                $
        &  $\frac{1}{4}(q+1)                                        $
        &  $\frac{1}{4}                                              $
        &  $\frac{1}{2}(q+1)                                $
        &  $\frac{1}{2}                                        $
        &  $1                                                 $   \\
 $\rb$  &  $-\frac{1}{2}(q^3-1)                                      $
        &  $\frac{1}{2}                                              $
        &  $\frac{1}{2}                                              $
        &  $m_{1}: \ -\frac{1}{4}(q^2-1)                             $
        &  $\frac{1}{4}(q+1)                                         $
        &  $\frac{1}{4}                                              $
        &  $-\frac{1}{2}(q-1)                                        $
        &  $\frac{1}{2}                                              $
        &                                                         \\
        &
        &
        &
        &  $m_{2}: \ -\frac{1}{4}                                    $
        &  $-\frac{1}{4}                                             $
        &  $\frac{1}{4}                                              $
        &
        &
        &                                                         \\
\hline
\end{tabular}
\end{center}
\vspace{0.25in}

%
\begin{center}
\begin{tabular}{||l|l|l|l|l|l|l|l|l|l||}  \hline
  Partition     & Class $H_{1}$ & Class $H_{2}$ & Class $H_{3}$
                & Class $I_{1}$ & Class $I_{2}$ & Class $J_{3}$
                & Class $K_{1}$ & Class $K_{2}$ & Class $L$       \\
                                       \hline

 $\rb$  &  $\frac{1}{6}(q+1)(q^2+q+1)                              $
        &  $\frac{1}{6}(2q+1)                                   $
        &  $\frac{1}{6}                                          $
        &  $\frac{1}{2}(q+1)                                $
        &  $\frac{1}{2}                                        $
        &  $1                                                         $
        &
        &
        &                                                            \\
 $\rc$  &  $-\frac{1}{2}(q^3-1)                              $
        &  $\frac{1}{2}                                   $
        &  $\frac{1}{2}                                $
        &  $-\frac{1}{2}(q-1)                                        $
        &  $\frac{1}{2}                                      $
        &
        &  $\frac{1}{2}(q^2+1)                               $
        &  $\frac{1}{2}                                      $
        &  $1$                                                \\
 $\re$  &  $\frac{1}{3}(q-1)(q^2-1)                             $
        &  $-\frac{1}{3}(q-1)                                   $
        &  $\frac{1}{3}                                $
        &
        &
        &
        &
        &
        &                                                  \\

 $\rf$  &
        &
        &
        &
        &
        &
        &$-\frac{1}{2}(q^2-1)                             $
        &$\frac{1}{2}                                     $
        &                                                  \\
                            \hline
\end{tabular}
\end{center}

\vspace{0.25in}
%
\begin{center}
\begin{tabular}{||l|l|l|l|l|l|l||}  \hline
  Partition     & Class $M_{1}$ & Class $M_{2}$ & Class $N$
                & Class $O$ & Class $P$ & Class $Q$             \\
                                       \hline
$\rd$   &  $\frac{1}{2}(q+1)                                       $
        &  $\frac{1}{2}                                   $
        &  $1                                          $
        &
        &
        &                                                            \\
$\re$   &  $-\frac{1}{2}(q-1)                                       $
        &  $\frac{1}{2}                                   $
        &
        &  $1$
        &
        &                                                            \\
$\rf$   &
        &
        &
        &
        &  $1$
        &                                                          \\
$\rf$   &
        &
        &
        &
        &
        &  $1$                                                        \\
                                       \hline
\end{tabular}
\end{center}
\end{landscape}

\begin{landscape}
%
\newcommand{\ra}{\mbox{$\rho_{1}$}}
\newcommand{\rb}{\mbox{$\rho_{2}$}}
\newcommand{\rc}{\mbox{$\rho_{3}$}}
\newcommand{\rd}{\mbox{$\rho_{4}$}}
\newcommand{\re}{\mbox{$\rho_{5}$}}
\newcommand{\rf}{\mbox{$\rho_{6}$}}
\newcommand{\rg}{\mbox{$\rho_{7}$}}
\newpage
\begin{center}
   Table 4. THE FUNCTIONS $\chi(m,e)$.
\end{center}
\begin{center}
\begin{tabular}{||c|c|c|c|c|c|c|c|c|c|c|c|c|c|c|c|c|c|c|c|c|c||}  \hline
  Partition  & $A_{1}^{'}$ & $A_{2}^{'}$ & $ A_{3}^{'}$ & $A_{4}^{'}$
             & $A_{5}^{'}$ & $A_{6}^{'}$ & $ A_{7}^{'}$
             & $B_{1}^{'}$ & $B_{2}^{'}$ & $ B_{3}^{'}$
             & $B_{4}^{'}$ & $B_{5}^{'}$
             & $C_{1}^{'}$ & $C_{2}^{'}$ & $ C_{3}^{'}$ & $C_{4}^{'}$
             & $C_{5}^{'}$ & $C_{6}^{'}$
             & $D_{1}^{'}$ & $D_{2}^{'}$ & $ D_{3}^{'}$      \\
                                    \hline
 $\ra$  &  $\frac{1}{120}        $
        &  $\frac{1}{30}         $
        &  $\frac{1}{24}         $
        &  $\frac{1}{20}         $
        &  $\frac{1}{24}         $
        &  $\frac{1}{30}         $
        &  $\frac{1}{120}        $
        &  $\frac{1}{24}         $
        &  $\frac{1}{8}          $
        &  $\frac{1}{12}         $
        &  $\frac{1}{8}          $
        &  $\frac{1}{24}         $
        &  $\frac{1}{12}         $
        &  $\frac{1}{6}          $
        &  $\frac{1}{12}         $
        &  $\frac{1}{12}         $
        &  $\frac{1}{6}        $
        &  $\frac{1}{12}        $
        &  $\frac{1}{6}         $
        &  $\frac{1}{3}         $
        &  $\frac{1}{6}         $
                          \\
\hline
\end{tabular}
\end{center}
\vspace{0.25in}

\begin{center}
\begin{tabular}{||c|c|c|c|c|c|c|c|c|c|c|c|c|c|c|c|c|c|c|c|c|c||}  \hline

  Partition
             & $E_{1}^{'}$ & $E_{2}^{'}$ & $ E_{3}^{'}$
             & $F_{1}^{'}$ & $F_{2}^{'}$
             & $G^{'}$ & $H_{1}^{'}$ & $ H_{2}^{'}$
             & $H_{3}^{'}$ & $I_{1}^{'}$
             & $I_{2}^{'}$ & $J^{'}$ & $ K_{1}^{'}$
             & $K_{2}^{'}$ & $L^{'}$ & $ M_{1}^{'}$
             & $M_{2}^{'}$ & $N^{'}$
             & $O^{'}$ & $P^{'}$ & $ Q^{'}$           \\
                                    \hline
 $\ra$
        &  $\frac{1}{4}         $
        &  $\frac{1}{4}        $
        &  $\frac{1}{4}         $
        &  $\frac{1}{4}          $
        &  $\frac{1}{2}         $
        &  $1                   $
        &  $0         $              &  $0$
        &  $0       $
        &  $0       $
        &  $ 0        $
        &  $ 0        $
        &  $ 0        $
        &  $ 0        $
        &  $ 0        $
        &  $ 0        $
        &  $ 0        $
        &  $ 0        $
        &  $ 0        $
        &  $ 0        $
        &  $ 0        $
          \\
 $\rb$

        &  $-\frac{1}{4}, -\frac{1}{4}         $
        &  $\frac{1}{4}, -\frac{1}{4}         $
        &  $\frac{1}{4}, \frac{1}{4}         $
        &  $-\frac{1}{2}        $
        &  $\frac{1}{2}         $
        &  $0                   $
        &  $\frac{1}{6}         $
        &  $\frac{1}{3}          $
        &  $\frac{1}{6}       $
        &  $\frac{1}{2}       $
        &  $\frac{1}{2}       $
        &  $1                 $
        &  $ 0        $
        &  $ 0        $
        &  $ 0        $
        &  $ 0        $
        &  $ 0        $
        &  $ 0        $
        &  $ 0        $
        &  $ 0        $
        &  $ 0        $
                                            \\
\hline
\end{tabular}
\end{center}
\end{landscape}
%
\newpage
\begin{center}
   Table 5. DUAL CLASSES (TYPES OF CHARACTERS)
\end{center}
\begin{center}
\begin{tabular}{||l|l|l|l|l||}  \hline
  Type      & Symbol & Range of Values & Number of & Degree of Characters \\
            &        & of $i,j,k,l,m$    & Characters  &                   \\
                                                       \hline
 $A_{7}^{'}$ & $g_{1i}^{5}$ & $1\leq i\leq (q-1)$ & $(q-1)$
          & $1$                                                         \\
$A_{6}^{'}$ & $g_{1i}^{4.1}$ & same & same
          & $(q+1)(q^2+1)$              \\
$A_{5}^{'}$ & $g_{1i}^{3.2}$ & same & same
          & $q^2(q^4+q^3+q^2+q+1)$         \\
$A_{4}^{'}$ & $g_{1i}^{3.1^{2}}$ & same & same
          & $q^{3}(q^2+1)(q^2+q+1)$            \\
$A_{3}^{'}$ & $g_{1i}^{2^2.1}$ & same & same
          & $q^{4}(q^4+q^3+q^2+q+1)$        \\
$A_{2}^{'}$ & $g_{1i}^{2.1^3}$ & same & same
          & $q^{6}(q^3+q^2+q+1)$     \\
\hline
\end{tabular}
\end{center}
\newpage
\begin{landscape}
\begin{center}
   Table 5. DUAL CLASSES (TYPES OF CHARACTERS) (Cont'd).
\end{center}
\begin{center}
\begin{tabular}{||l|l|l|l|l||}  \hline
  Type      & Symbol & Range of Values & Number of & Degree of Characters \\
            &        & of $i,j,k,l,m$    & Characters  &         \\
                                         \hline
$D_{1}^{'}$ & $g_{1i}^{1^3}.g_{1j}^{1}.g_{1k}^{1}$ & $1\leq i,j,k \leq (q-1)$
          & $\frac{1}{2}(q-1)(q-2)(q-3)$ &
          $q^{3}(q+1)(q^2+1)(q^4+q^3+q^2+q+1)$                \\
        &         & $i,j,k$ all &      &                            \\
        &         & distinct and &      &                            \\
        &         & $(i,j,k)=(i,k,j)$ &      &                            \\
$D_{2}^{'}$ & $g_{1i}^{2.1}.g_{1j}^{1}.g_{1k}^{1}$ & same
          & same &
          $q(q+1)^{2}(q^2+1)(q^4+q^3+q^2+q^+1)$                \\
$D_{3}^{'}$ & $g_{1i}^{3}.g_{1j}^{1}.g_{1k}^{1}$ & same
          & same &
          $q(q+1)(q^2+1)(q^4+q^3+q^2+q+1)$                \\
\hline
\end{tabular}
\end{center}
\end{landscape}
\newpage
%
\begin{landscape}
\begin{center}
   Table 5. DUAL CLASSES (TYPES OF CHARACTERS) (Cont'd).
\end{center}
\begin{center}
\begin{tabular}{||l|l|l|l|l||}  \hline
  Type      & Symbol & Range of Values & Number of & Degree of Characters \\
            &        & of $i,j,k,l,m$    & Characters &        \\
                                           \hline
$H_{1}^{'}$ & $g_{1i}^{1^3}.g_{2j}^{1} $ &
              $1\leq i \leq (q-1)$, $1\leq j\leq(q^2-2) $
          & $\frac{1}{2}q(q-1)^2 $ &
          $q^{3}(q^2+1)(q^5-1)$                \\
        &         & $j\not\equiv 0 \bmod (q+1) $ all distinct &     &                       \\
        &         & $j\sim jq \bmod(q^2-1) $  &    &                   \\
$H_{2}^{'}$ & $g_{1i}^{2.1}.g_{2j}^{1} $ &
               same
          &  same   &
          $q(q+1)(q^2+1)(q^5-1)$                \\
$H_{3}^{'}$ & $g_{1i}^{3}.g_{2j}^{1} $ &
               same
          &  same   &
          $(q^2+1)(q^5-1)$                \\
$I_{1}^{'}$ & $g_{1i}^{1^2}.g_{1j}^{1}.g_{2k}^{1} $ &
              $1\leq i,j \leq (q-1), i\neq j $
          & $\frac{1}{2}q(q-1)^{2}(q-2) $ &
          $q(q^2+1)(q^2+q+1)(q^5-1)$                \\
\hline
\end{tabular}
\end{center}
\end{landscape}

\newpage
%
\begin{landscape}
\begin{center}
   Table 5. DUAL CLASSES (TYPES OF CHARACTERS) (Cont'd).
\end{center}
\begin{center}
\begin{tabular}{||l|l|l|l|l||}  \hline
  Type  & Symbol & Range of Values & Number of & Number of Elements in \\
        &        & of $i,j,k,l,m$  & Characters & Degree of Characters \\
                                                       \hline
$M_{1}^{'}$ & $g_{1i}^{1^2}.g_{3j}^{1} $ &
              $1\leq i \leq (q-1), 1\leq j \leq (q^3-2) $   &
            $\frac{1}{3}q(q-1)^{2}(q+1) $ &
          $q(q^4-1)(q^5-1)$                \\
        &         & $ j\not\equiv 0 \bmod (q^2+q+1) $ &     &                       \\
        &         & $j\sim jq \sim jq^2 \bmod(q^3-1) $  &    &                   \\
$M_{2}^{'}$ & $g_{1i}^{2}.g_{3j}^{1} $ &
                same               &
                same               &
          $ (q^4-1)(q^5-1)$                \\
%
$N^{'}$ & $g_{1i}^{1}.g_{1j}^{1}.g_{3k}^{1} $ &
              $1\leq i,j \leq (q-1), i\neq j   $   &
            $\frac{1}{6}q(q-1)^{2}(q+1)(q-2) $ &
          $(q^4-1)(q^5-1)$                \\
        &         & $(i,j)=(j,i), 1\leq k\leq (q^3-2) $ &   &      \\
\hline
\end{tabular}
\end{center}
\end{landscape}
%
\newpage
Note: The contribution from the basic character corresponding to the missing
partitions is zero. The classes for which there is no contribution for any
of the partitions have been omitted.
\begin{center}
\begin{tabular}{||c|c||}               \hline           \hline
 \multicolumn{2}{||c||}
 {Table \# $6$}  \\
 \multicolumn{2}{||c||}
 {Basic Characters of Type $A^{'}$ }    \\    \hline
  Partition &  Contribution \\                     \hline
\multicolumn{2}{||c||}
  {Class A} \\                           \hline
      $\rho_{1}$   &   $120\epsilon_{1}^{5ia}$  \\
      $\rho_{2}$   &   $12\epsilon_{1}^{5ia}$  \\
      $\rho_{3}$   &   $8\epsilon_{1}^{5ia}$  \\
      $\rho_{4}$   &   $6\epsilon_{1}^{5ia}$  \\
      $\rho_{5}$   &   $6\epsilon_{1}^{5ia}$  \\
      $\rho_{6}$   &   $4\epsilon_{1}^{5ia}$  \\
      $\rho_{7}$   &   $5\epsilon_{1}^{5ia}$  \\              \hline
      \multicolumn{2}{||c||}{Class B}  \\                         \hline
      $\rho_{1}$   &   $120\epsilon_{1}^{i(4a+b)}$  \\
      $\rho_{2}$   &   $12\epsilon_{1}^{i(4a+b)}$  \\
      $\rho_{3}$   &   $8\epsilon_{1}^{i(4a+b)}$  \\
      $\rho_{4}$   &   $6\epsilon_{1}^{i(4a+b)}$  \\
      $\rho_{6}$   &   $4\epsilon_{1}^{i(4a+b)}$  \\          \hline
      \multicolumn{2}{||c||}{Class C}  \\                         \hline
      $\rho_{1}$   &   $120\epsilon_{1}^{i(3a+2b)}$  \\
      $\rho_{2}$   &   $m_{1}:12\epsilon_{1}^{i(3a+2b)}$  \\
                   &   $m_{2}:12\epsilon_{1}^{i(3a+2b)}$  \\
      $\rho_{3}$   &   $8\epsilon_{1}^{i(3a+2b)}$  \\
      $\rho_{4}$   &   $6\epsilon_{1}^{i(3a+2b)}$  \\
      $\rho_{5}$   &   $6\epsilon_{1}^{i(3a+2b)}$  \\         \hline
      \multicolumn{2}{||c||}{Class D}  \\                         \hline
      $\rho_{1}$   &   $120\epsilon_{1}^{i(3a+b+c)}$  \\
      $\rho_{2}$   &   $12\epsilon_{1}^{i(3a+b+c)}$  \\
      $\rho_{4}$   &   $120\epsilon_{1}^{i(3a+b+c)}$  \\     \hline
\end{tabular}
\end{center}
\newpage
\begin{center}
\begin{tabular}{||c|c||}               \hline
 \multicolumn{2}{||c||}
 {Table \# $6$}  \\
 \multicolumn{2}{||c||}
 {Basic Characters of Type $A^{'}$ (cont'd)}    \\    \hline
 Partition &  Contribution \\          \hline
\multicolumn{2}{||c||}
 {Class E} \\                           \hline         \hline
      $\rho_{1}$   &   $120\epsilon_{1}^{i(2a+2b+c)}$  \\
      $\rho_{2}$   &   $m_{1}:12\epsilon_{1}^{i(2a+2b+c)}$  \\
                   &   $m_{2}:12\epsilon_{1}^{i(2a+2b+c)}$  \\    \hline
      \multicolumn{2}{||c||}{Class F}  \\
\end{tabular}
\end{center}
\newpage
\begin{center}
\begin{tabular}{||c|c||}               \hline
 \multicolumn{2}{||c||}
 {Table \# $6$}  \\
 \multicolumn{2}{||c||}
 {Basic Characters of Type $A^{'}$ (cont'd)}    \\    \hline
 Partition &  Contribution \\                     \hline
 \multicolumn{2}{||c||}
 {Class N} \\                        \hline         \hline
      $\rho_{4}$   &   $6\epsilon_{1}^{i(a+b+c)}$  \\             \hline
      \multicolumn{2}{||c||}{Class O}  \\                             \hline
      $\rho_{5}$   &   $6\epsilon_{1}^{i(a+b)}$  \\               \hline
      \multicolumn{2}{||c||}{Class P}  \\                             \hline
      $\rho_{6}$   &   $4\epsilon_{1}^{i(a+b)}$  \\               \hline
      \multicolumn{2}{||c||}{Class Q}  \\                             \hline
      $\rho_{7}$   &   $5\epsilon_{1}^{ia}$  \\          \hline    \hline
\end{tabular}
\end{center}

%
\newpage
\begin{center}
\begin{tabular}{||c|l||}               \hline     \hline
 \multicolumn{2}{||c||}
 {Table \# $6$}  \\
 \multicolumn{2}{||c||}
 {Basic Characters of Type $B^{'}$ }    \\    \hline
  Partition &  \multicolumn{1}{c||} {Contribution} \\             \hline
\multicolumn{2}{||c||}
  {Class A} \\                           \hline
      $\rho_{1}$   &   $120\ea^{(4i+j)a}$  \\
      $\rho_{2}$   &   $12\ea^{(4i+j)a}$  \\
      $\rho_{3}$   &   $8\ea^{(4i+j)a}$  \\
      $\rho_{4}$   &   $6\ea^{(4i+j)a}$  \\
\end{tabular}
\end{center}
\newpage
\begin{center}
\begin{tabular}{||c|l||}               \hline     \hline
 \multicolumn{2}{||c||}
 {Table \# $6$}  \\
 \multicolumn{2}{||c||}
 {Basic Characters of Type $B^{'}$ (cont'd)}    \\    \hline
  Partition &  \multicolumn{1}{c||} {Contribution} \\             \hline
\multicolumn{2}{||c||}
  {Class F} \\                           \hline
      $\rho_{1}$   &   $ 12\ea^{2ia} \sum_{b,c,d} \ea^{i(b+c)+jd} + 48 \ea^{i(a+b+c+d)+ja}  $\\
%
   $\rho_{2}$  &  $2\ea^{2ia} \sum_{b,c,d} \ea^{i(b+c)jd}$ \\
        \hline
      \multicolumn{2}{||c||}{Class G}  \\                        \hline
  $\rho_{1}$  &  $\sum_{a,b,c,d} \ea^{i(a+b+c+d)+je}  $  \\
\end{tabular}
\end{center}
\newpage
\begin{center}
\begin{tabular}{||c|l||}               \hline     \hline
 \multicolumn{2}{||c||}
 {Table \# $6$}  \\
 \multicolumn{2}{||c||}
 {Basic Characters of Type $C^{'}$ }    \\    \hline
  Partition &  \multicolumn{1}{c||} {Contribution} \\             \hline
\multicolumn{2}{||c||}
  {Class A} \\                           \hline
      $\rho_{1}$   &   $120\ea^{(3i+2j)a}$  \\
      $\rho_{2}$   &   $m_{1}: 12\ea^{(3i+2j)a}$  \\
                   &   $m_{2}: 12\ea^{(3i+2j)a}$ \\
      $\rho_{3}$   &   $8\ea^{(3i+2j)a}$  \\
\hline
\end{tabular}
\end{center}
\newpage
\begin{center}
\begin{tabular}{||c|l||}               \hline     \hline
 \multicolumn{2}{||c||}
 {Table \# $6$}  \\
 \multicolumn{2}{||c||}
 {Basic Characters of Type $C^{'}$ (cont'd)}    \\    \hline
  Partition &  \multicolumn{1}{c||} {Contribution} \\             \hline
\multicolumn{2}{||c||}
  {Class E} \\                           \hline
      $\rho_{1}$   &   $24\ea^{i(2a+b)+j(a+c)}+
                        12\ea^{i(2a+c)+2jb}+$     \\
                   &   $24\ea^{i(a+2b)+j(a+c)}  +
                         4\ea^{i(a+b+c)+j(a+b)} +
                           \ea^{i(2b+c)+2ja}       $ \\
   $\rho_{2}$  &  $m_{1}:$       (i) $12\ea^{i(2a+c)+2jb}$ \\
               & \hspace{0.25in} (ii) $12\ea^{i(2b+c)+2ja}  $ \\
\end{tabular}
\end{center}
\newpage
\begin{center}
\begin{tabular}{||c|l||}               \hline     \hline
 \multicolumn{2}{||c||}
 {Table \# $6$}  \\
 \multicolumn{2}{||c||}
 {Basic Characters of Type $C^{'}$ (cont'd)}    \\    \hline
  Partition &  \multicolumn{1}{c||} {Contribution} \\         \hline
      \multicolumn{2}{||c||}{Class K}  \\                        \hline
      $\rho_{3}$   &   $8\ea^{i(a+b)+jb}$  \\           \hline
\multicolumn{2}{||c||}
  {Class L} \\                           \hline
      $\rho_{3}$   &   $4\ea^{i(a+b)+jc} +
                        4\ea^{i(a+c)+jb}   $ \\         \hline
%
%
      \multicolumn{2}{||c||}{Class M}  \\                        \hline
      $\rho_{4}$   &   $6\ea^{ib+2ja}$  \\
      $\rho_{5}$   &   $6\ea^{ib+2ja}$  \\              \hline
      \multicolumn{2}{||c||}{Class N}  \\                        \hline
      $\rho_{4}$   &   $6\ea^{ic+j(a+b)}$  \\           \hline
%
%
      \multicolumn{2}{||c||}{Class O}  \\                        \hline
      $\rho_{5}$   &   $6\ea^{ib+ja}$  \\    \hline  \hline
%
\end{tabular}
\end{center}
%
%
\newpage
\begin{center}
\begin{tabular}{||c|l||}               \hline     \hline
 \multicolumn{2}{||c||}
 {Table \# $6$}  \\
 \multicolumn{2}{||c||}
 {Basic Characters of Type $D^{'}$ }    \\    \hline
  Partition &  \multicolumn{1}{c||} {Contribution} \\             \hline
\multicolumn{2}{||c||}
  {Class A} \\                           \hline
      $\rho_{1}$   &   $120\ea^{(3i+j+k)a}$  \\
      $\rho_{2}$   &   $12\ea^{(3i+j+k)a}$  \\
      $\rho_{3}$   &   $6\ea^{(3i+j+k)a}$  \\             \hline
      \multicolumn{2}{||c||}{Class B}  \\
\hline
\end{tabular}
\end{center}
%
\newpage
\begin{center}
\begin{tabular}{||c|l||}               \hline     \hline
 \multicolumn{2}{||c||}
 {Table \# $6$}  \\
 \multicolumn{2}{||c||}
 {Basic Characters of Type $D^{'}$ (cont'd)}    \\    \hline
  Partition &  \multicolumn{1}{c||} {Contribution} \\             \hline
\multicolumn{2}{||c||}
  {Class G} \\                           \hline
      $\rho_{1}$   &   $\sum_{a,b,c,d,e}\ea^{i(a+b+c)+jd+ke} $ \\     \hline
      \multicolumn{2}{||c||}{Class H}  \\                       \hline
      $\rho_{2}$   &  $12\ea^{i(a+b)+ja+ka}$  \\                     \hline
      \multicolumn{2}{||c||}{Class I}  \\     \hline
      $\rho_{2}$   &  $2\ea^{ic}\sum_{i,j,k}\ea^{ib+(j+k)a}  $ \\    \hline
      \multicolumn{2}{||c||}{Class J}  \\                 \hline
      $\rho_{2}$   &  $2\ea^{id}\sum_{i,j,k}\ea^{ia+jb+kc}  $ \\     \hline
%
      \multicolumn{2}{||c||}{Class M}  \\                   \hline
      $\rho_{4}$   &  $6\ea^{ib+ja+ka}  $ \\                         \hline
%
      \multicolumn{2}{||c||}{Class N}  \\                   \hline
      $\rho_{4}$   &  $3\ea^{ic+ja+kb} + 3\ea^{ic+jb+ka}  $ \\ \hline \hline
%
\end{tabular}
\end{center}
%
%
\newpage
\begin{center}
\begin{tabular}{||c|l||}               \hline     \hline
 \multicolumn{2}{||c||}
 {Table \# $6$}  \\
 \multicolumn{2}{||c||}
 {Basic Characters of Type $E^{'}$ }    \\    \hline
  Partition &  \multicolumn{1}{c||} {Contribution} \\             \hline
\multicolumn{2}{||c||}
  {Class A} \\                           \hline
      $\rho_{1}$   &   $120\ea^{(2i+2j+k)a}$  \\
      $\rho_{2}$   & $m_{1}:  12\ea^{(2i+2j+k)a}$  \\
                   & $m_{2}:  12\ea^{(2i+2j+k)a} $        \\
      $\rho_{3}$   &   $8\ea^{(2i+2j+k)a}$  \\
\hline
\end{tabular}%
\end{center}
\newpage
\begin{center}
{\bf REFERENCES}
\end{center}

1. J. A. Green, ``The characters of general linear groups." {\em Trans. Amer.
                 Math. Soc.} 80(1955), 402-447.

2. S. K. Gupta, Ph.D. thesis, Case Institute of Technology 1971.

3. D. Littlewood, {\em The theory of groups characters and matrix
                   representation of groups.} Oxford. 1958.

4. D. Littlewood, ``On certain symmetric functions." $Proc. London Math.
                   Soc.$(3) 11 (1961), 485-498.

5. I.G. Macdonald, {\em Symmetric functions and Hall polynomials.}
                    Clarendon Press. Oxford 1979.

6. Alun O. Morriss, ``The characters of the groups $GL(n,q)$ ".
                     {\em Math. Z.} 81(1963), 733-742.

7. Alun O. Morriss, The multiplication of Hall Functions."
                     {\em Proc. London Math. Soc.} 3(13) 1963, 733-742.

8. W. Specht,     ``Die Characters der Symmetrischen Gruppe."
                   {\em Math. Z.} 73(1960) 312-329.

9. T. Springer,    {\em Characters of special groups.} Lecture Notes in
                    Mathematics. No. 131, Springer Verlag, 1970.

10. R. Steinberg,   Ph.D. Thesis, University of Toronto, 1948.

11. R. Steinberg,   ``The representations of $GL(3,q), GL(4,q), PGL(3,q),
                    PGL(4,q)$." {\em Canad. Jour. of Math.} 3 (1951),225-235.

12. R. Steinberg,   ``A geometric approach to the representations of the full
                    linear group over a Galois field." {\em Trans. Amer.
                    Math. Soc.} 71(1951), 274-282.

\end{document}